\documentclass[reqno,11pt]{amsart}

\usepackage{hyperref}
\usepackage{a4wide,amsthm,amssymb,mathrsfs,setspace,pstricks,booktabs,mathtools,amsmath,amsfonts}
\usepackage{graphicx}
\usepackage{xcolor}
\usepackage{array}
\newtheorem{theorem}{Theorem}[section]
\newtheorem{lemma}[theorem]{Lemma}
\newtheorem{proposition}[theorem]{Proposition}

\theoremstyle{definition}

\newtheorem{conjecture}[theorem]{Conjecture}

\theoremstyle{remark}
\newtheorem{remark}[theorem]{Remark}

\numberwithin{equation}{section}

\newcommand{\R}{\mathbb{R}}
\newcommand{\C}{\mathbb{C}}
\newcommand{\Rt}{\mathbb{R}^2}
\newcommand{\Rn}{\mathbb{R}^n}

\DeclareMathOperator{\spanop}{span}

\newcommand{\ip}[2]{\langle #1,#2\rangle}
\newcommand{\norm}[1]{\lVert #1\rVert}

\begin{document}
	
	\title [Planar Payne-P\'olya-Weinberger conjecture]{On the planar Payne-P\'olya-Weinberger conjecture for ratios of consecutive Dirichlet eigenvalues}
	
	\author{Long-Han Huang}
	\address{Department of Mathematical Sciences, Tsinghua University, Beijing 100084, China}
	\email{huanglh22@mails.tsinghua.edu.cn}

	\author{Zhen Song}    
	\address{Department of Mathematics, Northeastern University, Shenyang 110004, China}
	\email{songzhen@mail.neu.edu.cn}
	
	\author{Yuxuan Zhou}
	\address{Department of Mathematical Sciences, Tsinghua University, Beijing 100084, China}
	\email{zhou-yx22@mails.tsinghua.edu.cn}
	
	\author{Wenming Zou}
	\address{Department of Mathematical Sciences, Tsinghua University, Beijing 100084, China}
	\email{zou-wm@mail.tsinghua.edu.cn}


	
	\begin{abstract}
		Let $\Omega\subset\Rt$ be a bounded open set and $\lambda_k(\Omega)$ be the $k$-th Dirichlet eigenvalue of $\Omega$, counted with multiplicity. In this paper, we prove that for any $m\geq 1$,
		$$
		\frac{\lambda_{m+1}(\Omega)}{\lambda_m(\Omega)}\leq \frac{j_{1,1}^2}{j_{0,1}^2}.
		$$
		Here, $j_{p,k}$ denotes the $k$-th positive zero of the Bessel function $J_p$. Furthermore, if $\Omega$ is connected and $m\geq 2$, then the equality can never be attained. For each $m\geq1$, the constant $j_{1,1}^2/j_{0,1}^2$ is optimal even among bounded smooth domains. Therefore, our results completely resolve the long-standing Payne-P\'olya-Weinberger conjecture in the planar case.
	\end{abstract}

	\maketitle
	
	\section{Introduction}\label{section1}
\renewcommand{\thethm}{\Alph{thm}}
Let $\Omega\subset\Rn$ be a bounded domain with $n\geq 2$. Consider the classical Dirichlet eigenvalue problem for the Laplace operator on $\Omega$
\begin{equation}\label{dirichlet problem}
	\begin{cases}
			-\Delta u=\lambda u \quad &\hbox{in}\ \Omega, \\
			 u=0 \quad &\hbox{on}\ \partial\Omega. \\
		\end{cases}
\end{equation}
It is well-known that the spectrum of Dirichlet Laplacian consists of isolated eigenvalues of finite multiplicity accumulating to $\infty$. Let $\{\lambda_k(\Omega)\}_{k=1}^\infty$ be the Dirichlet eigenvalues of $\Omega$, ordered so that 
$$
0<\lambda_1(\Omega)<\lambda_2(\Omega)\leq\lambda_3(\Omega)\leq\cdots,\quad \lambda_k(\Omega)\to\infty\quad\text{as}\quad k\to\infty,
$$
with a corresponding orthogonal set of eigenfunctions denoted by $\{u_k\}_{k=1}^\infty\subset H_0^1(\Omega)$.

In 1956, Payne, P\'olya and Weinberger \cite{PaynePolyaWeinberger1956JMP} initiated the study of estimating the ratios of consecutive Dirichlet eigenvalues. More concretely, they proved in the planar case that the ratio of two consecutive Dirichlet eigenvalues satisfies
\begin{equation}\label{ppw inequality first}
	\frac{\lambda_{m+1}(\Omega)}{\lambda_m(\Omega)}\leq 3,\quad\forall\, m\geq 1,\quad \forall\, \Omega\subset \R^2.
\end{equation}
They also conjectured in \cite{PaynePolyaWeinberger1956JMP} that if both $m$ and the shape of domain are allowed to vary, then $\lambda_{m+1}/\lambda_{m}$ attains its maximum when $m=1$ and the domain is a disk. In 1969, Thompson \cite{Thompson1969SAM} proved that for any bounded domain $\Omega\subset\Rn$, $n\geq 2$, there holds
\begin{equation}\label{ppw inequality in dimension n by thompson}
	\frac{\lambda_2(\Omega)}{\lambda_1(\Omega)}\leq 1+\frac4n,
\end{equation}
thereby generalizing inequality \eqref{ppw inequality first} to dimension $n$ when $m=1$. Therefore, the general Payne-P\'olya-Weinberger Conjecture can be stated as follows:
\begin{conjecture}[Payne-P\'olya-Weinberger conjecture] \label{ppw conjecture general}
	Let $\Omega\subset\R^n$ be a bounded domain with $n\geq2$. Then
	\begin{equation}\label{ppw conjecture equation general}
			\frac{\lambda_{m+1}(\Omega)}{\lambda_m(\Omega)}\leq\left(\frac{j_{n / 2,1}}{j_{n / 2-1,1}}\right)^2,\quad \text{for }m=1,2,3,\cdots.
		\end{equation}
	Here, $j_{p,k}$ denotes the $k$-th positive zero of the Bessel function $J_p$.
\end{conjecture}
Due to the importance of this conjecture, many mathematicians have studied it extensively and obtained many important results, see \cite{Brands1964ARMA,Chiti1981BUMIA,Chiti1983SIAM,deVries1967ZNA,HileProtter1980Indiana,PaynePolyaWeinberger1955CR} for details. The first major breakthrough on this conjecture was due to Ashbaugh and Benguria, who proved that this conjecture holds for all $n\geq2$ in the case $m=1$ by combining Weinberger's center of mass result \cite{Weinberger1956JRMA}, Chiti's comparison result \cite{Chiti1982ZAMP}, and some fine estimates of Bessel zeros. We refer the interested reader to \cite{AshbaughBenguria1991BullAMS,AshbaughBenguria1992AnnMath,AshbaughBenguria1992CMP} for the detailed proof. By combining the method used in \cite{AshbaughBenguria1992AnnMath} and an ingenious nodal domain argument, Ashbaugh and Benguria proved that for any $n\geq2$, Conjecture \ref{ppw conjecture general} holds for $m=2$ and $m=3$ in \cite{AshbaughBenguria1991Duke} and \cite{AshbaughBenguria1993Edinburgh}, respectively. Also, the equality can never be attained if we require $\Omega$ to be a domain. Note that in \cite{AshbaughBenguria1993Edinburgh}, they actually proved that for any bounded domain $\Omega$ and any index $m\geq 2$ and $k\geq2$, there holds
\begin{equation}\label{ashbaugh benguria edin general}
	\frac{\lambda_{kl(m)}(\Omega)}{\lambda_m(\Omega)}\leq\sup_{\Omega'\subset\R^n}\frac{\lambda_{k}(\Omega')}{\lambda_1(\Omega')},
\end{equation}
where $l(m)$ denotes the number of nodal domains of the $m$-th Dirichlet eigenfunction $u_m$. Since $l(2)=2$ for every bounded domain $\Omega$, one has
\begin{equation}\label{ashbaugh benguria edin}
	\frac{\lambda_{4}(\Omega)}{\lambda_2(\Omega)}\leq\sup_{\Omega'\subset\R^n}\frac{\lambda_{2}(\Omega')}{\lambda_1(\Omega')}=\left(\frac{j_{n / 2,1}}{j_{n / 2-1,1}}\right)^2,
\end{equation}
and hence Conjecture \ref{ppw conjecture general} holds for $m=2$ and $m=3$. Note that for general domains, it seems impossible to give a good lower bound for $l(3)$ or $l(4)$, so the nodal domain method used in \cite{AshbaughBenguria1993Edinburgh} does not work for $m\geq4$. Very recently, by combining a quadratic-form identity with a Brouwer degree argument, the second author and Zhang \cite{SongZhang2026} proved that the Payne-P\'olya-Weinberger conjecture holds for any $n\geq 2$ in the case $m=4$. To our best knowledge, for $m\geq 5$, nothing better is known. We also refer the reader to \cite{AshbaughBenguria1993SIAM,Brands1964ARMA,ChenZheng2011JDE,HileProtter1980Indiana,LiTangZhang2026,Marcellini1980JDE,WangXia2021AnalPDE} and the references therein for recent developments concerning the stronger Payne-P\'olya-Weinberger summation conjecture and the Ashbaugh-Benguria reciprocal-gap conjecture, both of which imply the Payne-P\'olya-Weinberger conjecture for $m=1$.

In this paper, we focus on the planar case $n=2$. The main result of this paper is as follows.
\begin{theorem}\label{main result}
	Let $\Omega\subset\Rt$ be a bounded open set. For any $m\geq 1$, there holds
	\begin{equation}\label{main result equation}
			\frac{\lambda_{m+1}(\Omega)}{\lambda_m(\Omega)}\leq \frac{j_{1,1}^2}{j_{0,1}^2}.
	\end{equation}
Furthermore, if $\Omega$ is connected and $m\geq 2$, then the equality can never be attained.
\end{theorem}
\begin{remark}
	Following the construction described by Ashbaugh \cite{AshbaughBenguriaLaugesenWeidl2009OWR}, it is standard to prove that the constant $j_{1,1}^2/j_{0,1}^2$ in \eqref{main result equation} is optimal for any $m\geq 1$. Notably, for any $m\geq 1$, this constant is optimal even among bounded smooth domains. Indeed, let $D$ be the union of $m$ unit disks with pairwise disjoint closures. Then we know
	$$
	\lambda_1(D)=\cdots=\lambda_m(D)=j_{0,1}^2,\quad \lambda_{m+1}(D)=j_{1,1}^2.
	$$
	According to the approximation result established by He, Tang and Zhang \cite[Lemma 4.1]{HeTangZhang2026}, we know there exist bounded smooth domains $D_k\subset\Rt$ such that $\lambda_j(D_k)\to \lambda_j(D)$ for every fixed $j\geq1$. Taking $j=m,m+1$ gives
	$$
	\lim_{k\to\infty}\frac{\lambda_{m+1}(D_k)}{\lambda_m(D_k)}=\frac{\lambda_{m+1}(D)}{\lambda_m(D)}=\frac{j_{1,1}^2}{j_{0,1}^2},
	$$
	thus proving that the constant $j_{1,1}^2/j_{0,1}^2$ in \eqref{main result equation} is optimal among bounded smooth domain for any $m\geq1$.
\end{remark}
Theorem \ref{main result} hence completely resolves the Payne-P\'olya-Weinberger conjecture in the planar case. As mentioned earlier, we know very little about the number of nodal domains of the $m$-th Dirichlet eigenfunction, so it seems impossible to prove this conjecture merely by improving the arguments used in \cite{AshbaughBenguria1991Duke,AshbaughBenguria1993Edinburgh}. On the other hand, for higher eigenvalues, the number of orthogonality conditions grows with $m$, whereas the desired upper bound for $\lambda_{m+1}/\lambda_m$ remains the same. An intuitive idea is to project the test function onto the orthogonal complement of the first $m$ eigenfunctions so as to ensure that these conditions are satisfied. However, this may increase the Rayleigh quotient. To overcome this difficulty, we will use the radial multiplier developed by Ashbaugh and Benguria to construct two test functions with a common center. This is quite different from the choice of center made in \cite{SongZhang2026}. A winding-number argument then allows us to choose the center and the coefficients such that both test functions satisfy all $m$ orthogonality conditions. Finally, by the angular Fourier decomposition, we may deduce the desired upper bound for $\lambda_{m+1}/\lambda_m$ via the variational principle. 

We believe that the methods used in this paper will also be helpful for proving the Payne-P\'olya-Weinberger conjecture for $n\geq3$. Some partial progress in this direction has been made and we will investigate this further in future work. 

The rest of this paper is organized as follows. In Section \ref{section2}, we construct two test functions with a common center, establish the orthogonal conditions and reduce the proof of Theorem \ref{main result} to a two-function inequality. In Section \ref{section3}, we derive some estimates related to the Bessel functions. Finally in Section \ref{section4}, by introducing the angular Fourier decomposition, we give the proof of the two-function inequality presented in Section \ref{section2} and hence prove Theorem \ref{main result}.

\section{Reduction to a two-function auxiliary inequality}\label{section2}
From now on, let $\alpha=j_{0,1}$, $\beta=j_{1,1}$ and $d=\beta^2-\alpha^2$ for simplicity. Let $\Omega\subset\Rt$ be a bounded open set. Up to a domain dilation (see \cite[Lemma 2.3]{SongZhang2026} for instance), we may assume without loss of generality that $\lambda_{m}(\Omega)=\alpha^2$. In this section, we will reduce the proof of Theorem \ref{main result} to establishing a two-function auxiliary inequality.

To begin with, let $G\subset \Rt$ be an open set and $f,h\in H^1(G)$. For $\tau \in\R$, define
$$
q_\tau^{G}(f,h)=\int_G\left(\nabla f\cdot\nabla h-\tau fh\right)\,dx,
\quad q_\tau^{G}(f)=q_\tau^{G}(f,f).
$$
For $f,h\in H_0^1(\Omega)$, let $\widetilde{f},\widetilde{h}\in H^1(\Rt)$ denote the usual zero extensions of $f,h$. Then it is easy to know that
$$
q_\tau^{\Rt}(\widetilde{f},\widetilde{h})=q_\tau^{\Omega}(f,h),\quad q_\tau^{\Rt}(\widetilde{f})=q_\tau^{\Omega}(f).
$$
We shall always identify $f\in H_0^1(\Omega)$ with its zero extension $\widetilde{f}\in H^1(\Rt)$ whenever no confusion arises. Let $E_m=\spanop\{\varphi_1,\cdots,\varphi_m\}$, where $\varphi_1,\cdots,\varphi_m$ are $L^2(\Omega)$-orthonormal Dirichlet eigenfunctions corresponding to the first $m$ Dirichlet eigenvalues. We then know $q_{\alpha^2}^\Omega$ is non-positive on $E_m$. Indeed, for $u=\sum_{j=1}^m a_j\varphi_j$, we have
$$
q_{\alpha^2}^\Omega(u)=\sum_{j=1}^m(\lambda_j(\Omega)-\alpha^2)a_j^2\leq 0.
$$
To estimate $\lambda_{m+1}(\Omega)$, recall that from the min-max principle, we know
\begin{equation}\label{min-max}
	\lambda_{m+1}(\Omega)=\inf_{\substack{0\neq z\in H_0^1(\Omega)\\ z\perp E_m}}\frac{\int_\Omega|\nabla z|^2\,dx}{\int_\Omega z^2\,dx}.
\end{equation}
Therefore, to prove \eqref{main result equation}, it suffices to find a non-zero function $z\in H_0^1(\Omega)$ such that $z$ is orthogonal to $E_m$ and 
$$
q_{\beta^2}^\Omega(z)=\int_\Omega\left(|\nabla z|^2-\beta^2 z^2\right)\,dx\leq 0.
$$
This is the starting point of this paper.

Consider the following radial Bessel quotient
\begin{equation}\label{g}
	g(r)=\begin{cases}J_1(\beta r)/J_0(\alpha r),&0\le r<1,\\ \gamma,&r\ge1,\end{cases}
\end{equation}
where
$$
\gamma=\lim_{r\uparrow 1}J_1(\beta r)/J_0(\alpha r).
$$
This function has already been extensively studied by Ashbaugh and Benguria in their seminal paper \cite{AshbaughBenguria1992AnnMath}. By the results in \cite{AshbaughBenguria1992AnnMath,LiTangZhang2026,SongZhang2026}, we know $g\in C^1([0,\infty))$, $g(0)=0$, $g(1)=\gamma>0$, and $g'(r)>0$ for $r\in(0,1)$. Moreover, the Bessel quotient $g$ satisfies the following asymptotic properties:
\begin{equation}
	g(r)=\frac{\beta}{2}r+O(r^3),\quad r\downarrow0,
	\label{g origin behavior}
\end{equation}
\begin{equation}
	g(r)=\gamma\left(1-\frac{d-1}{6}(1-r)^2
	+O((1-r)^3)\right),\quad r\uparrow1,
	\label{g boundary behavior}
\end{equation}
and
\begin{equation}
	g'(r)=\frac{\gamma(d-1)}3(1-r)+O((1-r)^2),\quad r\uparrow1.
	\label{gprime boundary behavior}
\end{equation}
Recall that by \cite[Lemma 3.5]{LiTangZhang2026} and \cite[Lemma 3.1]{SongZhang2026}, we know $d>1$. Define
$$
B(r)\coloneqq\frac{g(r)^2}{\gamma^2},
$$
then by the above estimates, we know $0<B(r)<1$ for $0<r<1$.

For a center $A\in\Rt$, write $r=|x-A|$. Let $\phi_A(x)$ be the first positive Dirichlet eigenfunction of $B_1(A)$ such that $\|\phi_A\|_{L^2(B_1(A))}=1$. Extend $\phi_A$ to $\Rt$ by zero. It is easy to know that
\begin{equation}\label{phi}
	\phi_A(x)=\phi(r)=\begin{cases}c_\phi J_0(\alpha r),&0\le r<1,\\0,&r\ge1,\end{cases}
\end{equation}
where $c_\phi$ is a positive constant that ensures $\|\phi_A\|_{L^2(B_1(A))}=1$. Define the vector field
\begin{equation}\label{F_A}
	F_A(x)\coloneqq \begin{cases}
		g(r)\dfrac{x-A}{r},&x\ne A,\\
		0,&x=A,
	\end{cases}
\end{equation}
and write $F_A=(F_{A,1},F_{A,2})$ explicitly. We know from \cite[Lemma 3.1]{SongZhang2026} that $F_A\in C^1(\Rt,\Rt)\cap W^{1,\infty}(\Rt,\Rt)$ and $\|F_A\|_{W^{1,\infty}(\Rt,\Rt)}$ is uniformly bounded. For $f\in H^1(\Rt)$, define
\begin{equation}\label{R_A}
	R_A(f)\coloneqq q_{\alpha^2}^{\R^2}(f)-\frac1{\gamma^2}\sum_{i=1}^2q_{\beta^2}^{\R^2}(F_{A,i}f).
\end{equation}
The following energy identity has been used to prove the Payne-P\'olya-Weinberger conjecture for $m=4$ in \cite{SongZhang2026}.
\begin{proposition}\label{energy identity proposition}
	Fix $A\in\Rt$ and write $r=|x-A|$. For every $f\in H^1(\Rt)$, there holds
	\begin{equation}\label{energy identity equation}
		R_A(f) = \int_{B_1(A)} (1-B(r))\phi_A(x)^2 \left|\nabla\left(\frac{f}{\phi_A}\right)(x)\right|^2\,dx + \int_{\R^2\setminus B_1(A)} \left(d-\frac{1}{r^2}\right)f(x)^2\,dx.
	\end{equation}
	In particular, $R_A(f)\geq0$ with equality if and only if there exists $c\in\R$ such that $f=c\phi_A$ almost everywhere in $\Rt$.
\end{proposition}

\begin{proof}
	The identity \eqref{energy identity equation} is an extension of the radial remainder estimate in \cite{LiTangZhang2026}, and has already been established by Song and Zhang \cite[Proposition 3.2]{SongZhang2026}. Clearly, both integrals in \eqref{energy identity equation} are finite. Notice that by the property of $B$, we know $1-B(r)$ is positive in $B_1(A)$. For $r\ge1$, we know $d-1/r^2\geq d-1>0$. Therefore, it is easy to see that $R_A(f)\geq0$ for any $f\in H^1(\Rt)$. If $R_A(f)=0$, then
	$$
	\int_{\R^2\setminus B_1(A)} \left(d-\frac{1}{r^2}\right)f(x)^2\,dx=0,
	$$
	from which we know that $f=0$ almost everywhere in $\R^2\setminus B_1(A)$. On the other hand, since $\phi_A$ is positive in $B_1(A)$ and $B_1(A)$ is connected, we know $f/\phi_A=c$ almost everywhere in $B_1(A)$ for some constant $c\in\R$. Since $\phi_A(x)=0$ for $x\in \R^2\setminus B_1(A)$, we know $f=c\phi_A$ almost everywhere in $\R^2$. The proof is complete.
\end{proof}

\subsection{Simultaneous orthogonality of two test functions via a common center}

For $u,v\in H_0^1(\Omega)$ and $A\in \Rt$, define the test functions
\begin{equation}\label{test functions z1 and z2}
	z_1\coloneqq F_{A,1}u-F_{A,2}v,\quad z_2\coloneqq F_{A,2}u+F_{A,1}v.
\end{equation}
By \cite[Lemma 2.1]{SongZhang2026}, we know $z_1,z_2\in H_0^1(\Omega)$.

Briefly, we want to find a center $A\in\Rt$, such that $z_1$ and $z_2$ satisfy some orthogonal conditions. This idea goes back to Weinberger \cite{Weinberger1956JRMA}, who used it to establish the isoperimetric inequality for Neumann eigenvalue. We refer the reader to \cite{AshbaughBenguria1992AnnMath} for the application of this idea to Dirichlet eigenvalue problem.

\begin{proposition}[Simultaneous orthogonality]\label{center proposition}
	Let $E\subset H_0^1(\Omega;\R)$ be a subspace of dimension $m\ge1$, and let $P_E$ be the orthogonal projection of $L^2(\Omega;\R)$ onto $E$. Then there exist a center $A_*\in\Rt$ and $u,v\in E$ such that the test functions $z_1$ and $z_2$ defined in \eqref{test functions z1 and z2} with $A=A_*$ satisfy
	\begin{equation}\label{orthogonal condition of z1 and z2}
		P_Ez_1=P_Ez_2=0,\quad \norm{z_1}_{L^2(\Omega)}^2+\norm{z_2}_{L^2(\Omega)}^2>0.
	\end{equation}
\end{proposition}
\begin{proof}
	Let $e_1,\cdots,e_m$ be an $L^2(\Omega)$-orthonormal basis of $E$. Identify $\Rt$ with $\C$ and write $F_{A,1}+iF_{A,2}$ for the vector field defined in \eqref{F_A}. Define
	$$
	M(A)_{jk}\coloneqq\int_\Omega (F_{A,1}+iF_{A,2})e_ke_j\,dx,\quad D(A)\coloneqq\det M(A).
	$$
	Hence $M(A)$ is a $m\times m$ matrix. For $u,v\in E$, we may write $u+iv=\sum_{k=1}^m c_k e_k$ for some $\mathbf c=(c_1,\ldots,c_m)^T\in\mathbb C^m$. Then, it is easy to verify that
	$$
	(M(A)\mathbf c)_j=\int_\Omega(z_1+iz_2)e_j\,dx.
	$$
	Therefore, we know $P_Ez_1=P_Ez_2=0$ if and only if $M(A)\mathbf c=0$. Write $M=M(A)=M_1+iM_2$, where $M_1$ and $M_2$ are both real. Furthermore, by the property of determinant, we know
	\begin{equation}\label{determinant of D(A)}
		\det\begin{pmatrix}M_1&-M_2\\ M_2&M_1\end{pmatrix}=|D(A)|^2.
	\end{equation}
	Note that both $M$ and $D$ depend continuously on $A$ since $F_0$ is globally Lipschitz continuous in $\Rt$ and $F_A=F_0(\cdot-A)$.
	
	Take $L>0$ such that $\Omega\subset B_L(0)$ and $A=R\omega$ for some $R>L+1$ and $\omega=e^{i\vartheta}$. Then for $x\in\Omega$, we know $r=|x-A|\geq R-L>1$ and $|r-R|\leq |x|$. Since $g(r)=\gamma$ for all $r\geq1$, we know
	$$
	h_A(x)\coloneqq F_{A,1}(x)+iF_{A,2}(x)+\gamma\omega=\gamma\left(\frac{x-R\omega}{r}+\omega\right)=\gamma\frac{x+(r-R)\omega}{r}.
	$$
	Furthermore, we have
	$$
	|h_A(x)|\leq \gamma\frac{|x|+|r-R|}{r}\leq \frac{2\gamma L}{R-L}\eqqcolon\eta_R,
	\quad \forall\, x\in\Omega.
	$$
	For an $m\times m$ complex matrix $N$, recall that the operator norm is defined by
	$$
	\|N\|_{\mathrm{op}}\coloneqq\sup_{\mathbf c\in\mathbb C^m\setminus\{0\}}
	\frac{\|N\mathbf c\|_2}{\|\mathbf c\|_2},
	\quad
	\|\mathbf c\|_2\coloneqq\left(\sum_{k=1}^m|c_k|^2\right)^{1/2}.
	$$
	
	Let $P$ be the orthogonal projection in $L^2(\Omega,\mathbb C)$	onto the complex span of $e_1,\cdots,e_m$. For $\mathbf c=(c_1,\cdots,c_m)^{T}$, take $f=\sum_{k=1}^m c_ke_k$. By the definition of $M(A)$, it is easy to find out that
	$$
	\int_\Omega h_A f e_j\,dx=(M(A)\mathbf c)_j+\gamma\omega c_j,\quad 1\leq j\leq m.
	$$
	Therefore, we know $(M(A)+\gamma\omega I_m)\mathbf c$ is exactly the coefficient vector of $P(h_Af)$ with respect to the orthonormal basis $e_1,\cdots,e_m$. It follows that
	$$
	\|(M(A)+\gamma\omega I_m)\mathbf c\|_2=\|P(h_Af)\|_{L^2(\Omega)}\leq \|h_Af\|_{L^2(\Omega)}\leq\eta_R\|f\|_{L^2(\Omega)}=\eta_R\|\mathbf c\|_2,
	$$
	from which we know
	$$
	\|M(A)+\gamma\omega I_m\|_{\mathrm{op}}\leq\eta_R.
	$$
	For $1\leq j\leq m$, let $\mathbf a_j$ be the $j$-th column of $M(A)$ and $\mathbf b_j$ be the $j$-th column of $-\gamma\omega I_m$. Then we know
	$$
	\|\mathbf a_j-\mathbf b_j\|_2\leq \eta_R
	$$
	by the definition of operator norm. On the other hand, since $|\omega|=1$, we know $\|\mathbf b_j\|_2=\gamma$ and thus
	$$
	\|\mathbf a_j\|_2\leq\|\mathbf b_j\|_2+\|\mathbf a_j-\mathbf b_j\|_2\leq\gamma+\eta_R
	$$
	For $0\leq k\leq m$, let $\Delta_k$ be the determinant obtained by replacing the first $k$ columns of $M(A)$ by the first $k$ columns of $-\gamma\omega I_m$. A direct calculation shows that
	$$
	\Delta_k=\det(\mathbf b_1,\ldots,\mathbf b_k,\mathbf a_{k+1},\ldots,\mathbf a_m),\quad\Delta_0=D(A),\quad \Delta_m=(-\gamma\omega)^m.
	$$
	Moreover, for $1\leq k\leq m$, we know by Hadamard's inequality that
	$$
	\begin{aligned}
		|\Delta_{k-1}-\Delta_k|&\leq\left(\prod_{j=1}^{k-1}\|\mathbf b_j\|_2\right)\|\mathbf a_k-\mathbf b_k\|_2
		\left(\prod_{j=k+1}^{m}\|\mathbf a_j\|_2\right)\\
		&\leq \gamma^{k-1}\eta_R(\gamma+\eta_R)^{m-k}\\
		&\leq \eta_R(\gamma+\eta_R)^{m-1}.
	\end{aligned}
	$$
	Since $A=R\omega=Re^{i\vartheta}$, we have
	$$
	\bigl|D(Re^{i\vartheta})-(-\gamma)^m e^{im\vartheta}\bigr|=|\Delta_0-\Delta_m|\leq\sum_{k=1}^m|\Delta_{k-1}-\Delta_k|\leq m(\gamma+\eta_R)^{m-1}\eta_R.
	$$
	Notice that for fixed $m$, $m(\gamma+\eta_R)^{m-1}\eta_R$ tends to $0$ as $R\to\infty$ uniformly in $\vartheta$. Thus by taking $R$ sufficiently large, we may assume that $m(\gamma+\eta_R)^{m-1}\eta_R\leq \gamma^m/2$ for any $\vartheta$.
	
	Let $b(\vartheta)=(-\gamma)^m e^{im\vartheta}$. Define
	$$
	H(s,\vartheta)\coloneqq b(\vartheta)+s\bigl(D(Re^{i\vartheta})-b(\vartheta)\bigr),\quad 0\leq s\leq1.
	$$
	It is easy to verify that $H(s,\vartheta)$ is a homotopy that lies in $\C\setminus\{0\}$ since
	$$
	|H(s,\vartheta)|\geq |b(\vartheta)|-s|D(Re^{i\vartheta})-b(\vartheta)|\geq \frac{\gamma^m}{2}.
	$$
	Therefore, the loop $\vartheta\mapsto D(Re^{i\vartheta})$ has the same winding number around zero as the loop $\vartheta\mapsto b(\vartheta)$, and the winding number equals $m$. If $D\neq0$ in $\overline{B_R(0)}$, consider the following homotopy
	$$
	(t,\vartheta)\mapsto D(tRe^{i\vartheta}),\quad 0\leq t\leq1.
	$$
	Then we know the loop $\vartheta\mapsto D(Re^{i\vartheta})$ has the same winding number around zero as the constant loop $\vartheta\mapsto D(0)$, and the winding number equals $0$, which is a contradiction. Since $D\neq0$ on $\partial B_R(0)$, we know there exists $A_*\in B_R(0)$ such that $D(A_*)=0$.
	
	Now, choose $0\neq\mathbf c=(c_1,\cdots,c_m)^T\in\mathbb C^m$ such that $M(A_*)\mathbf c=0$ and take $u+iv=\sum_{k=1}^m c_ke_k$. Then we know $P_Ez_1=P_Ez_2=0$. Moreover, since
	$$
	\int_\Omega (u^2+v^2)\,dx=\sum_{k=1}^m|c_k|^2>0,
	$$
	and $g(|x-A_*|)>0$ except at $x=A_*$, we know from the identity $z_1^2+z_2^2=g(|x-A_*|)^2(u^2+v^2)$ that
	$$
	\|z_1\|_{L^2(\Omega)}^2+\|z_2\|_{L^2(\Omega)}^2=\int_\Omega g(|x-A_*|)^2(u^2+v^2)\,dx>0.
	$$
	The proof is complete.
\end{proof}

In \cite{SongZhang2026}, Song and Zhang used the theory of Brouwer degree to impose constraints on positive and negative parts of the test functions separately, thereby obtaining two centers. The difference is that, here, $z_1$ and $z_2$ share one common center.

\subsection{The two-function inequality}

Now, we compute $q_{\beta^2}^{\Rt}(z_1)+q_{\beta^2}^{\Rt}(z_2)$, where $z_1$ and $z_2$ are given in \eqref{test functions z1 and z2}. For $u,v\in H^1(\R^2)$, let $(r,\theta)$ be the standard polar coordinates centered at $A$. Define
\begin{equation}\label{mixed term}
	\mathcal{J}_A(u,v)\coloneqq2\int_{\R^2}\frac{B(r)}{r^2} \bigl(u\partial_\theta v-v\partial_\theta u\bigr)\,dx.
\end{equation}
Since $g(0)=0$ and $g$ is Lipschitz continuous, we know $B(r)/r$ is bounded. As $|\partial_\theta v|\le r|\nabla v|$, we know from the Cauchy-Schwarz inequality that the integral in \eqref{mixed term} converges absolutely and the functional $\mathcal{J}_A$ is continuous in $H^1(\R^2)\times H^1(\R^2)$.

Firstly, suppose that $u$ and $v$ are smooth. Write $z_1+iz_2=g(r)e^{i\theta}(u+iv)$. By a direct computation, we know
$$
\partial_r(z_1+iz_2)=e^{i\theta}\bigl(\partial_r(gu)+i\partial_r(gv)\bigr)
$$
and hence
$$
|\partial_r z_1|^2+|\partial_r z_2|^2=|\partial_r(gu)|^2+|\partial_r(gv)|^2
$$
Recall that $F_{A,1}=g\cos\theta$ and $F_{A,2}=g\sin\theta$. Thus
$$
\sum_{i=1}^2|\partial_r(F_{A,i}u)|^2=|\partial_r(gu)|^2,\quad \sum_{i=1}^2|\partial_r(F_{A,i}v)|^2=|\partial_r(gv)|^2
$$
Now we consider the angular derivative. Since
$$
\partial_\theta\!\left(e^{i\theta}(u+iv)\right)=e^{i\theta}\bigl[(\partial_\theta u-v)+i(\partial_\theta v+u)\bigr],
$$
it is easy to show that
\begin{align*}
	|\partial_\theta z_1|^2+|\partial_\theta z_2|^2
	&=g(r)^2\bigl[(\partial_\theta u-v)^2
	+(\partial_\theta v+u)^2\bigr]\\
	&=g(r)^2\bigl[(\partial_\theta u)^2+u^2
	+(\partial_\theta v)^2+v^2
	+2(u\partial_\theta v-v\partial_\theta u)\bigr].
\end{align*}
Besides, by the definition of $F_A$, we have
$$
\sum_{i=1}^2|\partial_\theta(F_{A,i}u)|^2=g(r)^2((\partial_\theta u)^2+u^2),\quad\sum_{i=1}^2|\partial_\theta(F_{A,i}v)|^2=g(r)^2((\partial_\theta v)^2+v^2).
$$
Since $g^2=\gamma^2B$ and 
$$
z_1^2+z_2^2=g(r)^2(u^2+v^2)=\sum_{i=1}^2\bigl((F_{A,i}u)^2+(F_{A,i}v)^2\bigr),
$$
we know by noticing \eqref{R_A} that
\begin{equation}\label{relation of energy}
	\frac{q_{\beta^2}^{\R^2}(z_1)+q_{\beta^2}^{\R^2}(z_2)}{\gamma^2}=q_{\alpha^2}^{\R^2}(u)+q_{\alpha^2}^{\R^2}(v)-R_A(u)-R_A(v)+\mathcal{J}_A(u,v).
\end{equation}
By using a density argument, we know the identity \eqref{relation of energy} holds for functions in $H^1(\Rt)$. Since $\mathcal{J}_A(v,u)=-\mathcal{J}_A(u,v)$, we know the sign of $\mathcal{J}_A$ is not fixed. According to \eqref{relation of energy}, we know an upper bound for $\mathcal{J}_A(u,v)$ gives an upper bound for $q_{\beta^2}^{\R^2}(z_1)+q_{\beta^2}^{\R^2}(z_2)$, and the following proposition gives the desired estimate when $v$ is orthogonal to $\phi_A$.
\begin{proposition}[Two-function inequality]\label{two function proposition}
	Let $A\in\R^2$, and let $u,v\in H^1(\R^2)$ be compactly supported. If $\ip{v}{\phi_A}_{L^2(\R^2)}=0$, then
	\begin{equation}\label{two function equation}
		R_A(u)+R_A(v)+q_{\alpha^2}^{\R^2}(v)-\mathcal{J}_A(u,v)\geq0,
	\end{equation}
	and the equality holds if and only if $v=0$ and $u=c\phi_A$ for some $c\in\R$.
\end{proposition}
The proof of Proposition \ref{two function proposition} will be given in Section \ref{section4}. Proposition \ref{two function proposition} provides a very useful estimate and we firstly use this proposition to prove Theorem \ref{main result}.
\begin{proof}[Proof of Theorem \ref{main result}]
	As we mentioned earlier in this section, $\lambda_m(\Omega)=\alpha^2$ and $E_m\subset H_0^1(\Omega)$ is the subspace spanned by the first $m$ Dirichlet eigenfunctions. By Proposition \ref{center proposition}, we know there exist a center $A_*\in\Rt$ and $u^0,v^0\in E_m$ such that the corresponding test functions $z_1^0$ and $z_2^0$ defined by \eqref{test functions z1 and z2} satisfy the orthogonal and non-trivial conditions \eqref{orthogonal condition of z1 and z2}. More concretely, both $z_1^0$ and $z_2^0$ are orthogonal to $E_m$ and at least one of them is non-zero.
	
	To apply Proposition \ref{two function proposition}, we need to ensure that $\phi_{A_*}$ and the zero extension of $v^0$ are orthogonal in $L^2(\Rt)$. Set
	$$
	a=\int_{\R^2}\bigl(\widetilde{u^0}+i\widetilde{v^0}\bigr)\phi_{A_*}\,dx,\quad
	\xi=\begin{cases}
		\overline a/|a|, & a\ne0,\\
		1, & a=0.
	\end{cases}
	$$
	Let $u,v\in E_m$ be defined by $u+iv=\xi(u^0+iv^0)$. Then
	$$
	\int_{\Rt}(\widetilde u+i\widetilde v)\phi_{A_*}\,dx=\xi a=|a|\in\R.
	$$
	By taking the imaginary part, we have
	$$
	\ip{\widetilde v}{\phi_{A_*}}_{L^2(\Rt)}=0.
	$$
	Now, let $z_1,z_2$ be the test functions induced by $u,v$ with the same center $A_*$. By \eqref{test functions z1 and z2}, we have
	$$
	z_1+iz_2=(F_{A_*,1}+iF_{A_*,2})(u+iv)=\xi(z_1^0+iz_2^0),
	$$
	from which we know $z_1$ and $z_2$ are real linear combinations of $z_1^0$ and $z_2^0$ and thus both $z_1$ and $z_2$ are orthogonal to $E_m$. Since $|\xi|=1$, we have $z_1^2+z_2^2=(z_1^0)^2+(z_2^0)^2$. In particular,
	$$
	\norm{z_1}_{L^2(\Omega)}^2+\norm{z_2}_{L^2(\Omega)}^2>0.
	$$
	
	Note that the zero extensions $\widetilde u,\widetilde v\in H^1(\R^2)$ have compact support in $\overline\Omega$. Since for $h\in H_0^1(\Omega)$, $q_\tau^{\R^2}(\widetilde h)=q_\tau^\Omega(h)$, we know by \eqref{relation of energy} that
	$$
	q_{\alpha^2}^\Omega(u)+2q_{\alpha^2}^\Omega(v)-\frac{q_{\beta^2}^\Omega(z_1)+q_{\beta^2}^\Omega(z_2)}{\gamma^2}
	=R_{A_*}(\widetilde u)+R_{A_*}(\widetilde v)+q_{\alpha^2}^{\R^2}(\widetilde v)-\mathcal{J}_{A_*}(\widetilde u,\widetilde v).
	$$
	Since $\widetilde v\perp\phi_{A_*}$, we know the right-hand side of the above equation is non-negative by Proposition \ref{two function proposition}. Recall that $q_{\alpha^2}^\Omega$ is non-positive on $E_m$, we finally get
	\begin{equation}\label{energy inequality final}
		\frac{q_{\beta^2}^\Omega(z_1)+q_{\beta^2}^\Omega(z_2)}{\gamma^2}\leq q_{\alpha^2}^\Omega(u)+2q_{\alpha^2}^\Omega(v)\leq0.
	\end{equation}
	Since $z_1,z_2$ are orthogonal to $E_m$, by the variational principle, we have
	\begin{equation}\label{energy inequality final 2}
		q_{\beta^2}^\Omega(z_1)+q_{\beta^2}^\Omega(z_2)\geq \bigl(\lambda_{m+1}(\Omega)-\beta^2\bigr)\bigl(\norm{z_1}_{L^2(\Omega)}^2+\norm{z_2}_{L^2(\Omega)}^2\bigr),
	\end{equation}
	and therefore, $\lambda_{m+1}(\Omega)\le\beta^2$. Since $\lambda_{m+1}/\lambda_{m}$ is invariant under domain dilation, we know \eqref{main result equation} is valid for any bounded open set $\Omega\subset\Rt$.
	
	Now, we prove that if $\Omega$ is connected and $m\geq 2$, then the equality in \eqref{main result equation} cannot be attained. Suppose to the contrary that $\lambda_{m+1}(\Omega)/\lambda_{m}(\Omega)=\beta^2/\alpha^2$ for some bounded connected domain $\Omega$ and $m\geq2$. We may further assume that $\lambda_m(\Omega)=\alpha^2$ and $\lambda_{m+1}(\Omega)=\beta^2$. Then by \eqref{energy inequality final 2}, we know both inequalities in \eqref{energy inequality final} are now equalities, which implies that
	$$
	R_{A_*}(\widetilde u)+R_{A_*}(\widetilde v)+q_{\alpha^2}^{\R^2}(\widetilde v)-\mathcal{J}_{A_*}(\widetilde u,\widetilde v)=0
	$$
	and
	$$
	q_{\alpha^2}^\Omega(u)+2q_{\alpha^2}^\Omega(v)=0.
	$$
	By Proposition \ref{two function proposition}, we know $\widetilde v=0$ and $\widetilde u=c\phi_{A_*}$ for some $c\in\R$. Note that at least one of $z_1$ and $z_2$ is non-zero, we know $c\neq0$ and $u$ is non-zero in $\Omega$. Moreover, $u$ has constant sign in $\Omega$. Combining the fact that $v=0$, we know $q_{\alpha^2}^\Omega(u)=0$.
	
	Since $u\in E_m$, we may write $u=\sum_{j=1}^m a_j \varphi_j$. Then we have
	$$
	0=q_{\alpha^2}^\Omega(u)=\sum_{j=1}^m\bigl(\lambda_j(\Omega)-\alpha^2\bigr)a_j^2.
	$$
	Since $\Omega$ is connected and $m\geq 2$, we know $\lambda_1(\Omega)<\lambda_2(\Omega)\leq\cdots\leq\lambda_{m}(\Omega)=\alpha^2$. Therefore, we must have $a_1=0$. However, this is impossible since
	$$
	a_1=\int_\Omega u\varphi_1\,dx
	$$
	and both $u$ and $\varphi_1$ are non-zero and have constant sign in $\Omega$. The proof is complete.
\end{proof}

\section{Bessel estimates}\label{section3}
From now on, our aim is to prove Proposition \ref{two function proposition}. In this section, we firstly establish some basic estimates related to the Bessel functions, which will play an important role in the proof of Proposition \ref{two function proposition} in Section \ref{section4}.

Recall that in \cite{SongZhang2026}, Song and Zhang obtained that
$$
	1-B(r)=\frac{d-1}{3}(1-r)^2+O\bigl((1-r)^3\bigr),\quad r\uparrow1.
$$
This estimate is insufficient to prove Proposition \ref{two function proposition} and we shall derive a lower bound for $1-B(r)$ for all $0\leq r\leq 1$.

From now on, we will use the change of variable $s=r^2$ without further mention and define
$$
	b(s)\coloneqq B(\sqrt{s}),\quad p(s)\coloneqq 1-b(s),\quad 0\leq s\leq1.
$$
Recently, Li, Tang and Zhang \cite[(3.24)]{LiTangZhang2026} proved that $b(s)\leq s e^{1-s}$, which gives the following lower bound
$$
p(s)\geq 1-s e^{1-s}\geq \frac12(1-s)^2,\quad 0<s<1.
$$
In Lemma \ref{estimate of c(s) and p(s)}, we will give a new polynomial lower bound
$$
	p(s)\geq\frac{13}{20}(1-s)^2+\frac{17}{60}(1-s)^3,\quad 0\le s\leq1,
$$
which improves on the quadratic bound $\frac12(1-s)^2$.

We first give the estimate for the Bessel zeros.
\begin{lemma}\label{bessel zeros estimate}
	There holds
	\begin{equation}\label{estimate alpha and beta}
		\frac{289}{50}<\alpha^2<\frac{579}{100},\quad\frac{367}{25}<\beta^2<\frac{1469}{100}.
	\end{equation}
	Besides, for $k\geq 2$, there holds
	\begin{equation}\label{estimate kth root}
		0<\frac{\alpha^2}{j_{0,k}^2}<\frac{\beta^2}{j_{1,k}^2}<1.
	\end{equation}
\end{lemma}

\begin{proof}
	For $\kappa=1,2$, denote $\Lambda_\kappa=j_{\kappa-1,1}^2$. Then $\Lambda_1=\alpha^2$ and $\Lambda_2=\beta^2$ are exactly the first Dirichlet eigenvalues of the unit balls in $\Rt$ and $\R^4$, respectively. Firstly, we use a similar argument as in \cite{AshbaughBenguria1992CMP} to derive upper bounds for $\alpha^2$ and $\beta^2$. Note that in $\R^{2\kappa}$, $\kappa=1,2$, we have for $u=u(r)$ that
	$$
	-\Delta u=-u''(r)-\frac{2\kappa-1}{r}u'(r).
	$$
	By taking $s=r^2$ and writing $h(s)=h(r^2)=u(r)$, we know
	$$
	u'(r)=2rh'(s)=2\sqrt{s}h'(s),\quad u''(r)=2h'(s)+4sh''(s).
	$$
	Thus, define
	$$
	(\mathcal L_\kappa h)(s)\coloneqq -4\bigl(sh''(s)+\kappa h'(s)\bigr)
	=-4s^{1-\kappa}\bigl(s^\kappa h'(s)\bigr)',\quad 0<s<1,
	$$
	and consider the eigenvalue problem $\mathcal L_\kappa h=\Lambda h$ together with the boundary condition $h(1)=0$. It is well-known that the function $$\psi_\kappa(s)\coloneqq s^{-(\kappa-1)/2}J_{\kappa-1}\bigl(\sqrt{\Lambda_\kappa s}\bigr),\quad 0\leq s\leq 1,$$ is an eigenfunction with eigenvalue $\Lambda_\kappa$. We may assume that $\psi_\kappa'(1)=-1$ up to multiplication by a positive constant. Take $h\in C^\infty([0,1])$ with $h(1)=0$. Recall that
	$$
	-4\bigl(s^\kappa \psi_\kappa'(s)\bigr)'=\Lambda_\kappa s^{\kappa-1}\psi_\kappa.
	$$
	Multiplying both sides of the above equation by $h^2/\psi_\kappa$ and then integrating over $(0,1)$ yields that
	$$
	4\int_0^1s^\kappa(h')^2\,ds-\Lambda_\kappa\int_0^1s^{\kappa-1}h^2\,ds=4\int_0^1s^\kappa\left(h'-\frac{\psi_\kappa'}{\psi_\kappa}h\right)^2ds.
	$$
	Therefore, for $h\not\equiv0$, we have
	\begin{equation}\label{estimate for Lambdak}
		\Lambda_\kappa\leq\frac{4\int_0^1s^\kappa(h')^2\,ds}{\int_0^1s^{\kappa-1}h^2\,ds}.
	\end{equation}
	For $\kappa=1$ and $\kappa=2$, by taking $h(s)=(1-s)(1-2s/5)$ and $h(s)=(1-s)(1-4s/5+s^2/5)$, respectively, we have $\alpha^2\le295/51<579/100$ and $\beta^2\le1248/85<1469/100$ after straightforward calculations.
	
	Now, we use a comparison principle to derive the lower bounds. Let $P$ be a polynomial satisfying that $P(1)=0$, $P>0$ on $[0,1)$ and $(\mathcal L_\kappa-t)P>0$ on $(0,1)$ for some $t>0$. We claim that $\Lambda_\kappa>t$. Indeed, since $\psi_\kappa$ is positive in $(0,1)$, we have
	\begin{equation}\label{comparison Lambdak and t}
		0<\int_0^1 s^{\kappa-1}\psi_\kappa(\mathcal L_\kappa-t)P\,ds=\int_0^1 s^{\kappa-1}P(\mathcal L_\kappa-t)\psi_\kappa\,ds=(\Lambda_\kappa-t)\int_0^1 s^{\kappa-1}P\psi_\kappa\,ds.
	\end{equation}
	Since $s^{\kappa-1}P\psi_\kappa$ is positive in $(0,1)$, we have $\Lambda_\kappa>t$. Next, we need to construct the polynomial $P$. Suppose that $P(s)=\sum_{j=1}^N \sigma_j(1-s)^j$ with $\sigma_1=1$. Set $\sigma_0=0$. Then for $0\leq j\leq N-2$, it is not hard to verify that the coefficient of $(1-s)^j$ in $(\mathcal L_\kappa-t)P$ is
	$$
	-4(j+1)(j+2)\sigma_{j+2}+4(j+1)(j+\kappa)\sigma_{j+1}-t\sigma_j.
	$$
	By taking
	\begin{equation}\label{lower coefficient}
		\sigma_{j+2}=\frac{j+\kappa}{j+2}\sigma_{j+1}
		-\frac{t}{4(j+1)(j+2)}\sigma_j,\quad 0\le j\le N-2,
	\end{equation}
	we have
	\begin{equation}\label{higher coefficient}
		(\mathcal L_\kappa-t)P=(1-s)^{N-1}\bigl[4N(N+\kappa-1)\sigma_N-t\sigma_{N-1}-t\sigma_N(1-s)\bigr]\,.
	\end{equation}
	If $t$ and $\sigma_N$ are positive, then we know $4N(N+\kappa-1)\sigma_N-t\sigma_{N-1}-t\sigma_N(1-s)$ is increasing in $s$, so to make sure that $(\mathcal L_\kappa-t)P>0$ on $(0,1)$, it suffices to verify that $4N(N+\kappa-1)\sigma_N-t\sigma_{N-1}-t\sigma_N>0$.
	For $\kappa=1$, take $t=289/50$ and $N=5$. By \eqref{lower coefficient}, we know
	$$
	\sigma_2=\frac12,\qquad \sigma_3=\frac{8-t}{24},\quad\sigma_4=\frac{6-t}{24}=\frac{11}{1200},\quad	\sigma_5=\frac{t^2-72t+384}{1920}.
	$$
	Besides, we know $t^2-72t+384=6/5+121/2500$ and hence $\sigma_5>1/1600$. All the coefficients $\sigma_1,\cdots,\sigma_5$ are positive and thus $P>0$ on $[0,1)$. Moreover, we also have
	$$
	100\sigma_5-t\sigma_4-t\sigma_5=(100-t)\sigma_5-t\sigma_4>\frac{94}{1600}-\frac{11}{200}=\frac3{800}>0,
	$$
	and we finally get that
	$$
	\alpha^2>\frac{289}{50}.
	$$
	For $\kappa=2$, take $t=367/25$ and $N=6$. By \eqref{lower coefficient}, we have $\sigma_2=1$, $\sigma_3=1-t/24$,
	$$
		\sigma_4=1-\frac{t}{16}=\frac{33}{400},\quad
		\sigma_5=1-\frac{3t}{40}+\frac{t^2}{1920},\quad
		\sigma_6=1-\frac{t}{12}+\frac{t^2}{960}.
	$$
	By a direct calculation, we may verify that $\sigma_1,\cdots,\sigma_6$ are all positive. Specifically, 
	$$
	0<\sigma_5=\frac{13489}{1200000}<\frac{9}{800},\quad \sigma_6=\frac{689}{600000}>\frac{1}{900}.
	$$
	Moreover,
	$$
		(168-t)\sigma_6-t\sigma_5>\frac{153}{900}-\frac{135}{800}=\frac1{800}>0.
	$$
	Therefore, we have
	$$
	\beta^2>\frac{367}{25}.
	$$
	
	Finally, the inequality \eqref{estimate kth root} has been established by Ashbaugh and Benguria, see \cite[Theorem 3.1]{AshbaughBenguria1992AnnMath} for instance. As stated by Ashbaugh and Benguria in \cite{AshbaughBenguria1992CMP}, we know that this result (and even a better one) was also obtained independently by Lorch. We refer the interested reader to \cite{Lorch1977,Lorch1980SIAM} for the proof.
\end{proof}

To proceed, we define the logarithmic derivatives for $0<s<1$ as follows:
\begin{equation}\label{a(s) and delta(s)}
	a(s)=-2\frac{d}{ds}\log\frac{J_0(\alpha\sqrt s)}{1-s},\quad
	\delta(s)=-2\frac{d}{ds}\log\frac{g(\sqrt s)}{\gamma\sqrt s}.
\end{equation}
The values of $a$ and $\delta$ at the endpoints $s=0$ and $s=1$ are defined by one-sided extensions. The following lemma provides several useful properties of $a$, $\delta$, and some other associated functions.

\begin{lemma}\label{lemma on delta and a}
	The function $\delta$ and all of its derivatives are strictly positive on $[0,1]$. Moreover,
	\begin{equation}\label{delta and a information}
		\delta(0)=\frac{\beta^2}{4}-\frac{\alpha^2}{2},\quad
		\delta(1)=1,\quad
		a(s)\ge a(0)=\frac{\alpha^2}{2}-2>\frac{89}{100}.
	\end{equation}
	Furthermore, we have
	\begin{equation}\label{bound of b}
		b(s)=s\exp\left(\int_s^1\delta(t)\,dt\right),\quad s\leq b(s)\leq se^{1-s}\leq1.
	\end{equation}
	For~$0<r<1$, there holds
	\begin{equation}\label{log derivative identity}
		\frac{rg'(r)}{g(r)}=1-s\delta(s),\quad
		-\frac{r\phi'(r)}{\phi(r)}=s\left(\frac2{1-s}+a(s)\right).
	\end{equation}
\end{lemma}
\begin{proof}
	Denote for $k\ge2$ that
	$$
		a_k=\frac{\alpha^2}{j_{0,k}^2},\quad
		b_k=\frac{\beta^2}{j_{1,k}^2}.
	$$
	By the product formula (see \cite{NISTDLMF10} for instance), we have
	$$
	\Pi_0(s)\coloneqq\frac{J_0(\alpha\sqrt s)}{1-s}=\prod_{k\ge2}(1-a_ks),\quad0<s<1,
	$$
	and
	$$
	\Pi_1(s)\coloneqq\frac{2J_1(\beta\sqrt s)}{\beta\sqrt s(1-s)}=\prod_{k\ge2}(1-b_ks),\quad 0<s<1.
	$$
	Since $0<a_k,b_k<1$ for $k\geq2$, both products are positive in $(0,1)$. By the expansion of Bessel function near $0$, we know $\Pi_0(0)=\Pi_1(0)=1$. Also, $\Pi_0$ and $\Pi_1$ can be extended smoothly to $s=1$ since $\alpha$ and $\beta$ are simple zeros of $J_0$ and $J_1$. For $0\leq x<1$, we have $x\leq -\log(1-x)$. By this elementary inequality, we know for $N\geq 2$ that
	$$
	\begin{aligned}
		\frac12\sum_{k=2}^N(a_k+b_k)&\leq -\sum_{k=2}^N\bigl[\log(1-a_k/2)+\log(1-b_k/2)\bigr]\\
		&\leq -\log\Pi_0(1/2)-\log\Pi_1(1/2).
	\end{aligned}
	$$
	Therefore, we know
	$$
	\sum_{k\ge2}(a_k+b_k)<\infty.
	$$
	By \eqref{estimate kth root}, we have $0<a_k<b_k\leq b_2<1$. For $x=a_k$ or $x=b_k$ and a non-negative integer $n$, we then know
	$$
	0<\frac{d^n}{ds^n}\frac{x}{1-xs}=n!\left(\frac{x}{1-xs}\right)^{n+1}\leq\frac{n!}{(1-b_2)^{n+1}}x,
	\quad 0\leq s\leq1.
	$$
	Consequently, for fixed $n\geq0$, the logarithmic derivative series
	$$
	\sum_{k\geq2}\frac{d^n}{ds^n}\log(1-a_ks),\quad\sum_{k\geq2}\frac{d^n}{ds^n}\log(1-b_ks)
	$$
	converge uniformly on $[0,1]$. Indeed, for $x=a_k$ or $x=b_k$ and $n\geq1$, we have
	$$ 
	\left|\frac{d^n}{ds^n}\log(1-xs)\right|=\left|-(n-1)!\left(\frac{x}{1-xs}\right)^n\right|\leq \frac{(n-1)!}{(1-b_2)^n}x,
	\quad 0\leq s\leq1,
	$$
	and for $x=a_k$ or $x=b_k$ and $n=0$, we have
	$$
	|\log(1-xs)|=-\log(1-xs)\leq \frac{xs}{1-xs}\leq \frac{x}{1-b_2}.
	$$
	Then the uniform convergence comes from the Weierstrass's Theorem.

	For the partial products
	$$
	P_{0,N}(s)\coloneqq\prod_{k=2}^N(1-a_ks),\quad P_{1,N}(s)\coloneqq\prod_{k=2}^N(1-b_ks),
	$$
	one has
	$$
	\log P_{0,N}(0)=\log P_{1,N}(0)=0.
	$$
	Since the logarithmic derivative series converge uniformly, we know for every $n\geq0$ and $i=0,1$,
	$$
	\frac{d^n}{ds^n}\log\Pi_i(s)=\lim_{N\to\infty}\frac{d^n}{ds^n}\log P_{i,N}(s)=\sum_{k\ge2}\frac{d^n}{ds^n}\log(1-c_{i,k}s),
	$$
	where $c_{0,k}=a_k$ and $c_{1,k}=b_k$. Hence the products $\Pi_0$ and $\Pi_1$ may be differentiated logarithmically term by term.

	By the definition of $g$, we know
	$$
	\frac{g(\sqrt s)}{\gamma\sqrt s}=\frac{\beta}{2\gamma}\,\frac{\Pi_1(s)}{\Pi_0(s)}.
	$$
	Then by \eqref{a(s) and delta(s)}, we know
	$$
	a(s)=-2\frac{\Pi_0'(s)}{\Pi_0(s)},\quad \delta(s)=-2\left(\frac{\Pi_1'(s)}{\Pi_1(s)}-\frac{\Pi_0'(s)}{\Pi_0(s)}\right).
	$$
	Direct substitution yields
	\begin{equation}\label{a(s) and delta(s) series}
		a(s)=2\sum_{k\ge2}\frac{a_k}{1-a_ks},\quad \delta(s)=2\sum_{k\ge2}\left(\frac{b_k}{1-b_ks}-\frac{a_k}{1-a_ks}\right),
	\end{equation}
	from which we know both $a$ and $\delta$ can be extended to the endpoints. The first identity also follows from Li, Tang and Zhang
	\cite[(3.10)]{LiTangZhang2026} and the second identity is equivalent to the expansion of $g'/g$ given by Ashbaugh and Benguria in \cite{AshbaughBenguria1992AnnMath}. For fixed $s\in[0,1]$, the map $x\mapsto x/(1-xs)$ is positive and strictly increasing in $(0,1)$. Since $b_k>a_k$, we know $\delta^{(n)}(s)>0$ for all $n\ge0$ and $s\in[0,1]$ since $\delta$ can be differentiated term by term by the above argument.
	
	At $s=0$, by the power series expansion of Bessel function (see \cite{NISTDLMF10} for instance), we have
	$$
	\Pi_0(s)=1+\left(1-\frac{\alpha^2}{4}\right)s+O(s^2),\quad
	\Pi_1(s)=1+\left(1-\frac{\beta^2}{8}\right)s+O(s^2).
	$$
	Since $\Pi_0'(0)=-\sum_{k\geq2}a_k$ and $\Pi_1'(0)=-\sum_{k\geq2}b_k$, we immediately get
	$$
		\sum_{k\geq2}a_k=\frac{\alpha^2}{4}-1,\quad\sum_{k\geq2}b_k=\frac{\beta^2}{8}-1.
	$$
	Therefore, by \eqref{a(s) and delta(s) series} and Lemma \ref{bessel zeros estimate}, we have
	$$
	\delta(0)=\frac{\beta^2}{4}-\frac{\alpha^2}{2},\quad a(s)\ge a(0)=\frac{\alpha^2}{2}-2>\frac{89}{100}.
	$$
	
	For $0<r<1$ and $s=r^2$, by the chain rule and the definition of $\phi(r)$, we have
	$$
	a(s)=-\frac{\phi'(r)}{r\phi(r)}-\frac2{1-s},\quad\delta(s)=\frac1s-\frac{g'(r)}{r\,g(r)}
	$$
	after a direct computation. This is exactly \eqref{log derivative identity}. Since $g(1)=\gamma>0$ and $g'(1)=0$, we obtain $\delta(1)=1$ by letting $r\uparrow 1$ in the above formula.
	
	Finally, 
	$$
	\frac{b(s)}{s}=\left(\frac{g(\sqrt{s})}{\gamma\sqrt{s}}\right)^2,
	$$
	so
	$$
	\left(\log\frac{b(s)}{s}\right)'=-\delta(s).
	$$
	Integrating from $s$ to $1$ yields that
	$$
		b(s)=s\exp\left(\int_s^1\delta(t)\,dt\right).
	$$
	Here, we have used the fact that $b(1)=1$. Since $\delta(s)$ is positive and increasing with $\delta(1)=1$, we know for $0<s\leq 1$ that
	$$
	s\leq b(s)\leq se^{1-s}\leq 1,
	$$
	where the last inequality follows from $\log s\leq s-1$. Finally, the continuity at $s=0$ gives \eqref{bound of b} on the whole interval $[0,1]$.
\end{proof}
By \eqref{delta and a information} and \eqref{bound of b}, we have $b(1)=1$ and $b'(1)=0$. To estimate $p(s)$, define
$$
C(s)\coloneqq\frac{b'(s)}{1-s},\quad H(s)\coloneqq\frac{1-\delta(s)}{1-s},\quad 0\leq s <1.
$$
Since $\delta(1)=1$, $H(s)$ can be extended continuously to $s=1$ with $H(1)=\delta'(1)$.

\begin{lemma}\label{estimate of c(s) and p(s)}
	The function $C$ can be extended continuously to $s=1$. Moreover, we have
	\begin{equation}\label{estimate of c(s) equation}
		C(s)\geq\frac{13}{10}+\frac{17}{20}(1-s),\quad 0\le s\leq1.
	\end{equation}
	Consequently,
	\begin{equation}\label{estimate of p(s) equation}
		p(s)\geq p_0(s)\coloneqq\frac{13}{20}(1-s)^2+\frac{17}{60}(1-s)^3,\quad 0\le s\leq1.
	\end{equation}
\end{lemma}

\begin{proof}
	Let $\psi_1,\psi_2$ be the normalized eigenfunctions given in the proof of Lemma \ref{bessel zeros estimate}. Obviously, the ratio $\psi_2/\psi_1$ is a constant multiple of $g(\sqrt{s})/\sqrt{s}$. Firstly, we need to determine the exact value of this constant. Since $\psi_1(1)=\psi_2(1)=0$ and $\psi_1'(1)=\psi_2'(1)=-1$, we know $\psi_2(s)/\psi_1(s)$ tends to $1$ as $s\uparrow 1$. Since $g(1)=\gamma$, we obtain that
	$$
	\frac{\psi_2(s)}{\psi_1(s)}=\frac{g(\sqrt s)}{\gamma\sqrt s},\quad 0<s<1.
	$$
	Notice that $s\psi_\kappa''+\kappa\psi_\kappa'+(\Lambda_\kappa/4)\psi_\kappa=0$ and $s=1$ is an ordinary point of this equation, so by a direct calculation, we may find out that
	$$
	\begin{aligned}
		\frac{\psi_1(s)}{1-s}&=1+\frac{1-s}{2}+\frac{8-\alpha^2}{24}(1-s)^2+\frac{6-\alpha^2}{24}(1-s)^3+O\bigl((1-s)^4\bigr),\\
		\frac{\psi_2(s)}{1-s}&=1+(1-s)+\frac{24-\beta^2}{24}(1-s)^2+\frac{16-\beta^2}{16}(1-s)^3+O\bigl((1-s)^4\bigr).
	\end{aligned}
	$$
	These two quotients are analytic and equal to $1$ at $s=1$. Using $d=\beta^2-\alpha^2$, we have
	$$
	\log\frac{\psi_2(s)}{\psi_1(s)}=\frac{1-s}{2}+\frac{7-d}{24}(1-s)^2+\frac{10-d}{48}(1-s)^3+O\bigl((1-s)^4\bigr).
	$$
	By \eqref{a(s) and delta(s)}, we know
	$$
	\delta(s)=-2\frac{d}{ds}\log\frac{\psi_2(s)}{\psi_1(s)}
	$$
	By taking derivative, we obtain
	$$
	c_1\coloneqq\delta'(1)=\frac{d-7}{6},\quad c_2\coloneqq\frac{\delta''(1)}2=\frac{10-d}{8}.
	$$
	By \eqref{estimate alpha and beta}, we know $889/100<d<891/100$ and therefore, $31/100<c_1<8/25$ and $0<c_2<7/50$.
	
	Recall that the Taylor's formula gives
	$$
	\delta(s)=1-c_1(1-s)+\int_s^1(t-s)\delta''(t)\,dt,\quad 0\leq s\leq1.
	$$
	By Lemma \ref{lemma on delta and a}, we know $\delta''\geq 0$ and $\delta'''\geq 0$ on $[0,1]$. Thus, $0\leq \delta''(t)\leq \delta''(1)=2c_2$ for $0\leq t\leq 1$, and
	$$
	1-c_1(1-s)\le\delta(s)\le1-c_1(1-s)+c_2(1-s)^2,\quad  0\leq s\leq 1.
	$$
	Therefore, we have
	$$
	H(s)\geq c_1-c_2(1-s)\ge\frac{31}{100}-\frac7{50}(1-s)
	$$
	for $0\leq s<1$ by the definition of $H(s)$. The inequality at $s=1$ follows from the continuity. By integrating from $s$ to $1$, we obtain that
	$$
	\int_s^1\delta(t)\,dt\geq(1-s)-\frac{c_1}{2}(1-s)^2 \geq(1-s)-\frac4{25}(1-s)^2.
	$$
	
	Differentiating \eqref{bound of b} gives
	$$
	b'(s)=\exp\left(\int_s^1\delta(t)\,dt\right)(1-s\delta(s)).
	$$
	Since $1-s\delta(s)=(1-s)(1+sH(s))$, we know
	$$
		C(s)=\exp\left(\int_s^1\delta(t)\,dt\right)(1+sH(s)).
	$$
	Thus we know $C(s)$ can be extended continuously to $s=1$ with $C(1)=1+c_1=(d-1)/6$. By combining the estimates for $H(s)$ and $\int_s^1\delta(t)\,dt$, we obtain that
	$$
	C(s)\geq\exp\left((1-s)-\frac4{25}(1-s)^2\right)\left[1+s\left(\frac{17}{100}+\frac7{50}s\right)\right].
	$$
	For $0\leq s\leq 1$, we have
	$$
	(1-s)-2\log\!\left(1+\frac{1-s}{2}\right)=\int_s^1\frac{1-t}{3-t}\,dt\geq\frac{(1-s)^2}{6}\ge\frac4{25}(1-s)^2.
	$$
	and hence
	$$
	\exp\left((1-s)-\frac4{25}(1-s)^2\right)\geq\left(1+\frac{1-s}{2}\right)^2.
	$$
	Direct expansion yields that
	\begin{align*}
		&\left(1+\frac{1-s}{2}\right)^2	\left[1+s\left(\frac{17}{100}+\frac7{50}s\right)\right]-\left(\frac{13}{10}+\frac{17}{20}(1-s)\right)\\
		&\quad\quad =\frac{2-s}{100}+\frac{(1-s)^2}{400}\bigl[7+11(1-s)+14(1-s)^2\bigr]>0,
	\end{align*}
	from which we know \eqref{estimate of c(s) equation} is valid. Finally, since $b'(t)=(1-t)C(t)$ and $b(1)=1$, we know
	$$
	p(s)=\int_s^1(1-t)C(t)\,dt\geq\frac{13}{20}(1-s)^2+\frac{17}{60}(1-s)^3,
	$$
	which is exactly \eqref{estimate of p(s) equation}. The proof is complete.
\end{proof}

\section{The proof of two-function inequality}\label{section4}

In this section, we shall prove Proposition \ref{two function proposition} by using an angular Fourier decomposition. As will be shown below, we will use three different methods to deal with this problem for different frequencies (the notion of frequency will be introduced later).

\subsection{Angular Fourier decomposition}\label{subsection 4.1}

By translations, we only need to prove Proposition \ref{two function proposition} for $A=0$. Write $B_1=B_1(0)$ and $r=|x|$. To decompose $R_0(f)$ into angular Fourier components, we use the following integral representation
\begin{equation}\label{integral representation}
	R_0(f)=\int_{\R^2}\left[(1-B(r))|\nabla f(x)|^2+V(r)f(x)^2\right]\,dx,
\end{equation}
where
\begin{equation}\label{definition of V}
	V(r)=
	\begin{cases}
		-\alpha^2(1-B(r))
		-B'(r)\dfrac{\phi'(r)}{\phi(r)}, & 0<r<1,\\[6pt]
		d-\dfrac{1}{r^2}, & r>1.
	\end{cases}
\end{equation}
This integral representation was proved by Song and Zhang in \cite[(3.24)]{SongZhang2026} for $f\in C_0^\infty(\Rt)$. To extend this representation to $H^1(\R^2)$, we first prove that $V$ is bounded. 

Near $r=0$, by the expansion of Bessel functions, we know $B'(r)=O(r)$ and $\phi'(r)/\phi(r)=O(r)$. Note that $\phi$ has a simple zero at $r=1$, so
$$
\frac{\phi'(r)}{\phi(r)}=-\frac{1}{1-r}+O(1),\quad r\uparrow1.
$$
On the other hand, near $r=1$, we have $g'(r)=O(1-r)$ and hence $B'(r)=2g(r)g'(r)/\gamma^2=O(1-r)$. Hence, we know the product $B'(r)\phi'(r)/\phi(r)$ is bounded as $r\uparrow 1$. For $r>1$, we immediately know that $d-r^{-2}$ is bounded. Since $1-B$ and $V$ are all bounded, the integral representation for $R_0$ given above is continuous with respect to the  strong convergence in $H^1(\Rt)$. By \eqref{R_A}, we know $R_0$ is also continuous in $H^1(\Rt)$. Therefore, by a density argument, we know \eqref{integral representation} holds for every $f\in H^1(\Rt)$.

Let $u,v\in H^1(\Rt)$ be two functions that satisfy the hypothesis of Proposition \ref{two function proposition}. For $w\in\{u,v\}$, let $w_0$ denote the angular mean of $w$ and let $w_\ell^c,w_\ell^s$ denote its cosine and sine coefficients at frequency $\ell\ge1$. Then for almost every $r>0$, we have
$$
w_0(r)=\frac{1}{2\pi}\int_0^{2\pi}w(r,\theta)\,d\theta
$$
and
$$
w_\ell^c(r)=\frac{1}{\pi}\int_0^{2\pi}w(r,\theta)\cos(\ell\theta)\,d\theta,
\quad
w_\ell^s(r)=\frac{1}{\pi}\int_0^{2\pi}w(r,\theta)\sin(\ell\theta)\,d\theta,
\quad \ell\geq1.
$$
For $\ell\geq1$, write $w_\ell(r,\theta)=w_\ell^c(r)\cos(\ell\theta)+w_\ell^s(r)\sin(\ell\theta)$, and for any non-negative integer $N$, define the partial summation of $w$ as follows:
$$
S_N w\coloneqq\sum_{\ell=0}^N w_\ell.
$$
Firstly, we compute the quadratic terms of \eqref{two function equation} for $S_N u$ and $S_N v$.

Fix an integer $\ell \geq 0$. Let $Y_0(\theta)=1$ and let $Y_\ell(\theta)$ be either $\cos(\ell\theta)$ or $\sin(\ell\theta)$ for $\ell\geq1$. Then
$$
\|Y_\ell\|_{L^2(0,2\pi)}^2=\int_0^{2\pi}Y_\ell(\theta)^2\,d\theta
=\begin{cases}
	2\pi,&\ell=0,\\
	\pi,&\ell\ge1.
\end{cases}
$$
For each fixed angular frequency $\ell$, the space of radial Fourier coefficients consists of real locally absolutely continuous functions $f$ on $(0,\infty)$ such that
$$
\|f\|_\ell^2\coloneqq\int_0^\infty\left[|f'(r)|^2+\left(1+\frac{\ell^2}{r^2}\right)|f(r)|^2\right]r\,dr<\infty.
$$
For such $f$, the function $h(r,\theta)=f(r)Y_\ell(\theta)$ belongs to $H^1(\R^2)$. Define the radial quadratic forms by
$$
\mathscr Q_\ell(f)\coloneqq\frac{q_{\alpha^2}^{\R^2}(h)}{\|Y_\ell\|_{L^2(0,2\pi)}^2},
\quad
\mathscr R_\ell(f)\coloneqq\frac{R_0(h)}{\|Y_\ell\|_{L^2(0,2\pi)}^2}.
$$
To compute these forms, observe that $\partial_r h=f'(r)Y_\ell(\theta)$ and $\partial_\theta h=f(r)Y_\ell'(\theta)$. Since $-Y_\ell''=\ell^2Y_\ell$, we have
$$
\int_0^{2\pi}Y_\ell'(\theta)^2\,d\theta=\ell^2\|Y_\ell\|_{L^2(0,2\pi)}^2.
$$
By using the polar coordinate, we have
$$
\int_0^{2\pi}h^2\,d\theta=\|Y_\ell\|_{L^2(0,2\pi)}^2 f(r)^2
$$
and
$$
\int_0^{2\pi}|\nabla h|^2\,d\theta=\|Y_\ell\|_{L^2(0,2\pi)}^2\left[f'(r)^2+\frac{\ell^2}{r^2}f(r)^2\right].
$$
Therefore,
\begin{equation}\label{definition of Q in sec 4}
	\mathscr Q_\ell(f)=\int_0^\infty\left[f'(r)^2+\left(\frac{\ell^2}{r^2}-\alpha^2\right)f(r)^2\right]r\,dr.
\end{equation}
Note that both $1-B(r)$ and~$V(r)$ are independent of the angular variable $\theta$. Then by a similar substitution and \eqref{integral representation}, we have
\begin{equation}\label{definition of R in sec 4}
	\mathscr R_\ell(f)=\int_0^\infty\left[(1-B(r))f'(r)^2+\left(V(r)+(1-B(r))\frac{\ell^2}{r^2}\right)f(r)^2\right]r\,dr.
\end{equation}
In particular, for any $\ell\geq 1$, both $\mathscr Q_\ell(f)$ and $\mathscr R_\ell(f)$ are independent of the choice of $Y_\ell(\theta)$. By the orthogonality of the trigonometric basis, we have
$$
R_0(S_Nu)=2\pi\,\mathscr R_0(u_0)+\pi\sum_{\ell=1}^N\bigl[\mathscr R_\ell(u_\ell^c)+\mathscr R_\ell(u_\ell^s)\bigr],
$$
$$
R_0(S_Nv)=2\pi\,\mathscr R_0(v_0)+\pi\sum_{\ell=1}^N\bigl[\mathscr R_\ell(v_\ell^c)+\mathscr R_\ell(v_\ell^s)\bigr],
$$
and
$$
q_{\alpha^2}^{\R^2}(S_Nv)=2\pi\,\mathscr Q_0(v_0)+\pi\sum_{\ell=1}^N\bigl[\mathscr Q_\ell(v_\ell^c)+\mathscr Q_\ell(v_\ell^s)\bigr].
$$

It remains to compute the mixed term $\mathcal{J}_0(S_Nu,S_Nv)$. By \eqref{mixed term}, we know
$$
\mathcal{J}_0(S_Nu,S_Nv)=2\int_0^\infty\frac{B(r)}{r}\int_0^{2\pi}\bigl(S_Nu\,\partial_\theta S_Nv-S_Nv\,\partial_\theta S_Nu\bigr)\,d\theta\,dr.
$$
At frequency $\ell$, taking derivative with respect to $\theta$ gives
$$
u_\ell\partial_\theta v_\ell-v_\ell\partial_\theta u_\ell=\ell\bigl(u_\ell^c v_\ell^s-u_\ell^s v_\ell^c\bigr).
$$
Note that the right-hand side is independent of $\theta$. Therefore, by the orthogonality of the trigonometric basis, we have
$$
\mathcal{J}_0(S_Nu,S_Nv)=4\pi\sum_{\ell=1}^N\ell\int_0^\infty\frac{B(r)}{r}\bigl(u_\ell^c(r)v_\ell^s(r)-u_\ell^s(r)v_\ell^c(r)\bigr)\,dr.
$$

An important observation is that the finite sums $S_Nu,S_Nv$ also satisfy the hypotheses of Proposition \ref{two function proposition}. Indeed, for $w\in\{u,v\}$,
\[
S_Nw(x)=\frac{1}{2\pi}\int_0^{2\pi}
\left(1+2\sum_{\ell=1}^N\cos(\ell\psi)\right)w(\rho_{-\psi}x)\,d\psi,
\]
where $\rho_{-\psi}$ denotes the rotation around origin through the angle $-\psi$. Since rotations preserve the $H^1(\Rt)$ norm and $1+2\sum_{\ell=1}^N\cos(\ell\psi)$ is integrable, we know $S_Nw\in H^1(\Rt)$. Moreover, $S_Nw$ is supported in any disk centered at the origin that contains the support of $w$. Note that the angular mean of $S_Nv$ is $v_0$. Since $\phi$ is radial and vanishes outside $B_1$, we have
$$
\langle S_Nv,\phi\rangle_{L^2(\R^2)}=2\pi\int_0^1 v_0(r)\phi(r)r\,dr=\langle v,\phi\rangle_{L^2(\R^2)}=0.
$$

For each $\ell\geq 0$, define the following frequency component
$$
D_\ell\coloneqq R_0(u_\ell)+R_0(v_\ell)+q_{\alpha^2}^{\R^2}(v_\ell)-\mathcal{J}_0(u_\ell,v_\ell).
$$
Based on the previous calculations, there holds
$$
R_0(S_Nu)+R_0(S_Nv)+q_{\alpha^2}^{\R^2}(S_Nv)-\mathcal{J}_0(S_Nu,S_Nv)=\sum_{\ell=0}^N D_\ell.
$$
We will prove that $D_\ell\ge0$ for every $\ell\geq0$.

For $\ell=0$, both $u_0$ and $v_0$ are independent of $\theta$, so $\mathcal{J}_0(u_0,v_0)=0$. Therefore, it suffices to show that 
\begin{equation}\label{inequality l=0}
	D_0=2\pi\bigl[\mathscr R_0(u_0)+\mathscr R_0(v_0)+\mathscr Q_0(v_0)\bigr]\geq0
\end{equation}
under the condition $\langle v_0,\phi\rangle_{L^2(\R^2)}=0$ established above.

For $1\leq\ell\leq N$, we may write
\begin{equation}\label{group of Dl}
	\begin{aligned}
		\frac{D_\ell}{\pi}&=\left[\mathscr R_\ell(u_\ell^c)+\mathscr R_\ell(v_\ell^s)+\mathscr Q_\ell(v_\ell^s)-4\ell\int_0^\infty\frac{B(r)}{r}\,u_\ell^c v_\ell^s\,dr\right]\\
		&\quad+\left[\mathscr R_\ell(u_\ell^s)+\mathscr R_\ell(v_\ell^c)+\mathscr Q_\ell(v_\ell^c)+4\ell\int_0^\infty\frac{B(r)}{r}\,u_\ell^s v_\ell^c\,dr\right].
	\end{aligned}
\end{equation}
Therefore, it suffices to prove that
\begin{equation}\label{inequality higher l}
	\mathscr R_\ell(f)+\mathscr R_\ell(k)+\mathscr Q_\ell(k)
	-4\ell\int_0^\infty\frac{B(r)}{r}f(r)k(r)\,dr\geq0.
\end{equation}
Note that taking $(f,k)=(u_\ell^c,v_\ell^s)$ gives the first bracketed term in \eqref{group of Dl}, and taking $(f,k)=(u_\ell^s,-v_\ell^c)$ gives the second bracketed term in \eqref{group of Dl}.

Once \eqref{inequality l=0} and \eqref{inequality higher l} are proved, we know \eqref{two function equation} holds for $(S_Nu,S_Nv)$ for every $N$. To pass to $(u,v)$, we need to prove that $S_Nw\to w$ in $H^1(\Rt)$ for any $w\in H^1(\Rt)$. We first assume that $w\in C_0^\infty(\Rt)$. By using the Parseval's identity on circle $\partial B_r$ to $w$, $\partial_r w$ and $\partial_\theta w$, we have
$$
\|w\|_{H^1(\Rt)}^2=2\pi\|w_0\|_0^2+\pi\sum_{\ell=1}^\infty\left(\|w_\ell^c\|_\ell^2+\|w_\ell^s\|_\ell^2\right).
$$
Similarly, we have
\begin{equation}\label{upper bound Sn}
	\|S_Nw\|_{H^1(\Rt)}\le\|w\|_{H^1(\Rt)},
\end{equation}
and
$$
\|S_Nw-w\|_{H^1(\R^2)}\to0\quad\text{as}\quad N\to\infty.
$$
Since $C_0^\infty(\Rt)$ is dense in $H^1(\Rt)$ and the upper bound \eqref{upper bound Sn} is uniform in $N$, we know $S_Nw\to w$ in $H^1(\Rt)$ for any $w\in H^1(\Rt)$. By continuity of $R_0$, $q_{\alpha^2}^{\Rt}$ and $\mathcal{J}_0$, we know the desired inequality holds for $(u,v)$.

\subsection{The case $\ell=0$}\label{subsection 4.2}

Since $\mathcal{J}_0(u_0,v_0)=0$ and $R_0(u_0)\geq0$ by Proposition \ref{energy identity proposition}, it remains to estimate $R_0(v_0)+q_{\alpha^2}^{\Rt}(v_0)$. Let $e_0=1/\sqrt\pi$ be a positive constant that will be used for normalization. Since
$$
J_0'(r)=-J_1(r),\quad J_1'(r)=J_0(r)-\frac{1}{r}J_1(r),
$$
we may verify that
$$
\int_0^1 rJ_0(\alpha r)^2\,dr=\frac{J_1(\alpha)^2}{2}, \quad \int_0^1 rJ_0(\alpha r)\,dr=\frac{J_1(\alpha)}{\alpha},
$$
from which we know the normalized constant $c_\phi$ given in \eqref{phi} is $c_\phi=1/(\sqrt\pi J_1(\alpha))$, and hence $\ip{\phi}{e_0}_{L^2(B_1)}=2/\alpha$. Using $\ip{v_0}{\phi}_{L^2(B_1)}=0$, we obtain that
$$
\left|\ip{v_0}{e_0}_{L^2(B_1)}\right|^2 =\left|\ip{v_0}{e_0-(2/\alpha)\phi}_{L^2(B_1)}\right|^2\leq\left(1-\frac4{\alpha^2}\right)\norm{v_0}_{L^2(B_1)}^2.
$$
Denote
$$
\overline v_0=\frac{1}{\pi}\int_{B_1}v_0\,dx.
$$
Since $\|e_0\|_{L^2(B_1)}=1$, we know $\overline v_0=\langle v_0,e_0\rangle_{L^2(B_1)}e_0$ and hence
$$
\|v_0-\overline v_0\|_{L^2(B_1)}^2=\|v_0\|_{L^2(B_1)}^2-|\langle v_0,e_0\rangle_{L^2(B_1)}|^2
\geq \frac{4}{\alpha^2}\|v_0\|_{L^2(B_1)}^2.
$$

Now we consider the function $v_0-\overline v_0$, which is radial and has integral zero over $B_1$. Since the first positive Neumann eigenvalue in the radially symmetric space of $B_1$ is $\beta^2$, by the variational principle,
$$
	\int_{B_1}|\nabla v_0|^2\,dx=\int_{B_1}|\nabla(v_0-\overline v_0)|^2\,dx\geq \beta^2\norm{v_0-\overline v_0}_{L^2(B_1)}^2\geq \frac{4\beta^2}{\alpha^2}\norm{v_0}_{L^2(B_1)}^2,
$$
which gives
\begin{equation}\label{q inside B1}
	q_{\alpha^2}^{B_1}(v_0|_{B_1})\geq \left(\frac{4\beta^2}{\alpha^2}-\alpha^2\right)\norm{v_0}_{L^2(B_1)}^2.
\end{equation}
Here, the coefficient $4\beta^2/\alpha^2-\alpha^2$ is positive by \eqref{estimate alpha and beta}. The contribution of $\Rt\setminus\overline{B_1}$ to $R_0(v_0)+q_{\alpha^2}^{\R^2}(v_0)$ is
$$
\int_{\Rt\setminus B_1}\left[|\nabla v_0|^2+ \left(\beta^2-2\alpha^2-r^{-2}\right)v_0^2\right]dx\geq0,
$$
since the coefficient of $v_0^2$ is positive by \eqref{estimate alpha and beta}. By noting that $1-B$ is non-negative in $B_1$, we finally get $R_0(v_0)+q_{\alpha^2}^{\R^2}(v_0)\geq0$, and the equality holds if and only if $v_0=0$.

\subsection{The case $\ell=1,2$}\label{subsection 4.3}

Fix $\ell\in\{1,2\}$. We begin by reducing \eqref{inequality higher l} to two scalar inequalities. Since
$$
-4\ell fk=-2\ell f^2-2\ell k^2+2\ell(f-k)^2,
$$
the left-hand side of equation \eqref{inequality higher l} is the sum of
\begin{equation}\label{scalar pair of inequality}
	\mathscr R_\ell(f)-2\ell\int_0^\infty\frac{B(r)}r f^2\,dr,\quad	\mathscr R_\ell(k)+\mathscr Q_\ell(k)-2\ell\int_0^\infty\frac{B(r)}r k^2\,dr,
\end{equation}
and the following non-negative term
\[
2\ell\int_0^\infty\frac{B(r)}r(f-k)^2\,dr.
\]
Therefore, it suffices to prove that both terms in \eqref{scalar pair of inequality} are non-negative.

First, we need to construct a suitable comparison function $h_\ell$. Inside $B_1$, we know $r^\ell$ satisfies the following equation:
$$
-(r^\ell)''-\frac1r(r^\ell)'+\frac{\ell^2}{r^2}r^\ell=0.
$$

We take $h_\ell(r)=r^\ell e^{1-r^2}$ in $B_1$, and then extend it to the exterior region by a function that matches both the value and the first derivative of $h_\ell$ at $r=1$. Namely, 
\begin{equation}\label{definition of h for comparison}
	h_\ell(r)=
	\begin{cases}
		r^\ell e^{1-r^2}, & 0<r\le1,\\
		e^{1-r}, & r>1,\ \ell=1,\\
		1, & r>1,\ \ell=2.
	\end{cases}
\end{equation}
Indeed, $h_\ell(1)=1$ and $h_\ell'(1^-)=h_\ell'(1^+)=\ell-2$. By a direct calculation, we have for $0<r<1$ that
$$
\frac{h_\ell'(r)}{h_\ell(r)}=\frac{\ell}{r}-2r,\quad-\frac{(rh_\ell'(r))'}{r h_\ell(r)}+\frac{\ell^2}{r^2}=4(\ell+1)-4r^2.
$$

For $r>0$, $r\neq1$, define
$$
\begin{aligned}
	W_\ell(r)&\coloneqq-\frac{\bigl(r(1-B(r))h_\ell'(r)\bigr)'}{r h_\ell(r)}+V(r)+\frac{\ell^2(1-B(r))-2\ell B(r)}{r^2},\\
	E_\ell(r)&\coloneqq-\frac{\bigl(rh_\ell'(r)\bigr)'}{r h_\ell(r)}+\frac{\ell^2}{r^2}-\alpha^2.
\end{aligned}
$$
For $f,k\in C_0^\infty(0,\infty)$, we claim that
\begin{equation}\label{comparison identity for W and E}
	\begin{aligned}
		&\mathscr R_\ell(f)-2\ell\int_0^\infty\frac{B(r)}r f^2\,dr\\
		&\quad\quad=\int_0^\infty(1-B(r))h_\ell^2\left|\left(\frac{f}{h_\ell}\right)'\right|^2r\,dr+\int_0^\infty W_\ell(r)f^2\,r\,dr,\\[4pt]
		&\mathscr R_\ell(k)+\mathscr Q_\ell(k)-2\ell\int_0^\infty\frac{B(r)}r k^2\,dr\\
		&\quad\quad=\int_0^\infty(2-B(r))h_\ell^2\left|\left(\frac{k}{h_\ell}\right)'\right|^2r\,dr+\int_0^\infty\bigl(W_\ell(r)+E_\ell(r)\bigr)k^2\,r\,dr.
	\end{aligned}
\end{equation}
Indeed, since
$$
(1-B)f'^2=(1-B)h_\ell^2\left|\left(\frac{f}{h_\ell}\right)'\right|^2+(1-B)(f^2)'\frac{h_\ell'}{h_\ell}-(1-B)f^2\frac{(h_\ell')^2}{h_\ell^2},
$$
integrating by parts on $(0,1)$ and $(1,\infty)$ respectively gives
$$
	\int_0^\infty(1-B)f'^2r\,dr=\int_0^\infty(1-B)h_\ell^2\left|\left(\frac{f}{h_\ell}\right)'\right|^2r\,dr-\int_0^\infty\frac{\bigl(r(1-B)h_\ell'\bigr)'}{r h_\ell}f^2r\,dr.
$$
Since $h_\ell$, $h_\ell'$ and $B$ are all continuous at $r=1$, the boundary terms cancel. Substituting this into the definition of $\mathscr R_\ell$ gives the first identity of \eqref{comparison identity for W and E}. Similarly, one has
$$
\mathscr Q_\ell(k)=\int_0^\infty h_\ell^2\left|\left(\frac{k}{h_\ell}\right)'\right|^2r\,dr+\int_0^\infty E_\ell(r)k^2r\,dr,
$$
from which we get the second identity of \eqref{comparison identity for W and E}.

Since $0\leq B\leq1$, we know both
$$
\int_0^\infty(1-B(r))h_\ell^2\left|\left(\frac{f}{h_\ell}\right)'\right|^2r\,dr
$$
and
$$
\int_0^\infty(2-B(r))h_\ell^2\left|\left(\frac{k}{h_\ell}\right)'\right|^2r\,dr
$$
are non-negative. We now prove that $W_\ell$ and $W_\ell+E_\ell$ are positive in $(0,1)$ and $(1,\infty)$.

For $0<r<1$, the product rule gives
$$
-\frac{\bigl(r(1-B)h_\ell'\bigr)'}{r h_\ell}=-(1-B)\frac{(rh_\ell')'}{r h_\ell}+B'\frac{h_\ell'}{h_\ell}.
$$
By the definition of $h_\ell$ and $V$, we have
\begin{equation}\label{E and W inside}
	\begin{aligned}
		E_\ell(r)&=4(\ell+1)-4r^2-\alpha^2,\\
		W_\ell(r)&=(1-B(r))E_\ell(r)-B'(r)\left(\frac{\phi'(r)}{\phi(r)}-\frac{\ell}{r}+2r\right)-\frac{2\ell B(r)}{r^2}.
	\end{aligned}
\end{equation}

\begin{lemma}\label{positivity inside}
	For $\ell=1,2$ and $0<r<1$, there holds
	\[
	W_\ell(r)>0\quad\text{and}\quad W_\ell(r)+E_\ell(r)>0.
	\]
\end{lemma}

\begin{proof}
	For $0<r<1$, we write $s=r^2$ as in Section \ref{section3}. Since $B(r)=g(r)^2/\gamma^2$,
	\eqref{log derivative identity} gives
	$$
	rB'(r)=2B(r)\frac{rg'(r)}{g(r)}=2b(s)\bigl(1-s\delta(s)\bigr),
	\quad	
	-\frac{r\phi'(r)}{\phi(r)}=s\left(\frac{2}{1-s}+a(s)\right).
	$$
	Substituting these identities into \eqref{E and W inside}, we obtain that
	$$
	\begin{aligned}
		W_\ell(r)&=p(s)E_\ell(r)+2b(s)\bigl(1-s\delta(s)\bigr)\left(\frac{2}{1-s}+a(s)+\frac{\ell}{s}-2\right)-\frac{2\ell b(s)}{s}\\
		&=p(s)E_\ell(r)+2b(s)\left[\bigl(1-s\delta(s)\bigr)\left(\frac{2}{1-s}+a(s)-2\right)-\ell\delta(s)\right].
	\end{aligned}
	$$
	Since
	$$
	\frac{1-s\delta(s)}{1-s}=\delta(s)+H(s)=1+sH(s),
	$$
	one has
	$$
	W_\ell(r)=p(s)E_\ell(r)+b(s)\left[2(2-\ell)\delta(s)+4H(s)+2\bigl(a(s)-2\bigr)\bigl(1-s\delta(s)\bigr)\right].
	$$
	
	We first estimate $4H(s)+2\bigl(a(s)-2\bigr)\bigl(1-s\delta(s)\bigr)$, since it is independent of the frequency $\ell$. By \eqref{estimate alpha and beta} and \eqref{delta and a information}, we know
	$$
	\frac{31}{40}<\delta(0)<\frac{79}{100},\quad a(s)>\frac78.
	$$
	Since $\delta$ is strictly increasing, i.e., $\delta(s)\geq \delta(0)$, we know from $\delta(1)=1$ and the convexity of $\delta$ that
	$$
	\delta(s)\le(1-s)\delta(0)+s,\quad	H(s)\ge 1-\delta(0)>\frac{21}{100}.
	$$
	Now, $1-s\delta(s)=(1-s)(1+sH(s))\geq0$ and $s(1-s)\leq1/4$, and therefore,
	$$
	\begin{aligned}
		4H(s)+2\bigl(a(s)-2\bigr)\bigl(1-s\delta(s)\bigr)&\geq 4H(s)-\frac94(1-s)\bigl(1+sH(s)\bigr)\\
		&=\left(4-\frac94s(1-s)\right)H(s)-\frac94(1-s)\\
		&\geq \left(4-\frac{9}{16}\right)\frac{21}{100}-\frac94(1-s)\\
		&\geq\frac7{10}-\frac94(1-s).
	\end{aligned}
	$$
	
	For $\ell=1$, the remaining term is $2\delta(s)\geq 2\delta(0)>31/20$. Since
	$$
	\frac{31}{20}+\frac7{10}-\frac94(1-s)=\frac94s,
	$$
	we have
	$$
	W_1(r)\geq p(s)E_1(r)+\frac94 b(s)s.
	$$
	By Lemma \ref{estimate of c(s) and p(s)}, we know $p(s)\geq0$. Hence if $E_1(r)\geq 0$, we know from $b(s)s>0$ for $0<s<1$ that $W_1(r)>0$ and $W_1(r)+E_1(r)>0$. Suppose that $E_1(r)<0$ for some $0<r<1$. Since $p(s)=1-b(s)$,
	$$
	W_1(r)+E_1(r)\geq 2E_1(r)+b(s)\left(\frac94s-E_1(r)\right).
	$$
	By \eqref{bound of b}, we know $b(s)\geq s$. Moreover, 
	$$
	E_1(r)=8-4s-\alpha^2>2-4s,
	$$
	where the last inequality follows from the fact that $\alpha^2<6$. Therefore,
	$$
	W_1(r)+E_1(r)\geq (2-s)E_1(r)+\frac94s^2>(2-s)(2-4s)+\frac94s^2=\left(2-\frac52s\right)^2.
	$$
	It follows that $W_1(r)+E_1(r)>0$. Since $E_1(r)<0$, we also have $W_1(r)>0$.
	
	For $\ell=2$, the term $2(2-\ell)\delta(s)$ vanishes, and
	$$
	E_2(r)=12-4s-\alpha^2>2+4(1-s)>0.
	$$
	Therefore, we may find out that
	$$
	\begin{aligned}
		W_2(r)&\geq p(s)\bigl[2+4(1-s)\bigr]+\bigl(1-p(s)\bigr)\left[\frac7{10}-\frac94(1-s)\right]\\
		&=\frac7{10}-\frac94(1-s)+p(s)\left[\frac{13}{10}+\frac{25}{4}(1-s)\right].
	\end{aligned}
	$$
	By \eqref{estimate of p(s) equation},
	$$
	W_2(r)\geq\frac7{10}-\frac94(1-s)+\frac34(1-s)^2+\frac{15}{4}(1-s)^3.
	$$
	Since
	$$
	(1-s)^3-\frac{1-s}{3}+\frac{2}{27}=\left(s-\frac23\right)^2\left(\frac53-s\right)\geq0,
	$$
	we obtain that
	$$
	W_2(r)\geq\frac{19}{45}-(1-s)+\frac34(1-s)^2=\frac34\left(s-\frac13\right)^2+\frac4{45}>0,
	$$
	from which we also know that $W_2(r)+E_2(r)>0$ since $E_2(r)>0$. The proof is complete.
\end{proof}

For $r>1$, we have $B(r)=1$ and $V(r)=d-r^{-2}$, and hence
$$
W_1(r)=d-\frac{3}{r^2},\quad W_1(r)+E_1(r)=\beta^2-2\alpha^2-1+\frac{1}{r}-\frac{2}{r^2}
$$
and
$$
W_2(r)=d-\frac{5}{r^2},\quad W_2(r)+E_2(r)=\beta^2-2\alpha^2-\frac{1}{r^2}.
$$
By \eqref{estimate alpha and beta}, we have $d>5$ and $\beta^2-2\alpha^2>2$, from which we immediately obtain that $W_1(r)$, $W_1(r)+E_1(r)$, $W_2(r)$ and $W_2(r)+E_2(r)$ are all positive in $(1,\infty)$.

Summing up the above, we have proved that for $f,k\in C_0^\infty(0,\infty)$, there holds
$$
\begin{aligned}
	&\mathscr R_\ell(f)+\mathscr R_\ell(k)+\mathscr Q_\ell(k)-4\ell\int_0^\infty\frac{B(r)}{r}f(r)k(r)\,dr\\
	&\quad\quad\geq \int_0^\infty\left[W_\ell(r)f(r)^2+\bigl(W_\ell(r)+E_\ell(r)\bigr)k(r)^2\right]r\,dr.
\end{aligned}
$$
We need to extend this estimate to general radial Fourier coefficients. Let $w$ be a radial Fourier coefficient with $\|w\|_\ell<\infty$. Choose $0\leq\eta_j\leq 1$ in $C_0^\infty(0,\infty)$ such that $\eta_j=1$ on $[1/j,j]$, $\operatorname{supp}\eta_j\subset[1/(2j),2j]$, and
$$
r|\eta_j'(r)|\leq M,\quad j\geq1,\quad r>0,
$$
where $M>0$ is a constant independent of $j$. Since $(w-\eta_jw)'=(1-\eta_j)w'-\eta_j'w$ and $\ell\geq1$, we know
$$
\|w-\eta_jw\|_\ell^2\leq 2(1+M^2)\int_{(0,1/j)\cup(j,\infty)}\left[|w'(r)|^2+\left(1+\frac{\ell^2}{r^2}\right)|w(r)|^2\right]r\,dr\to0
$$
as $j\to\infty$. Therefore, for general radial Fourier coefficients $f$ and $k$, by using a standard mollifying argument, we may find $f_j,k_j\in C_0^\infty(0,\infty)$, such that
$$
\|f_j-f\|_\ell+\|k_j-k\|_\ell\rightarrow0.
$$
Passing to a subsequence, we may assume that $f_j\to f$ and $k_j\to k$ almost everywhere. By Fatou's lemma, we finally have
$$
\begin{aligned}
	&\int_0^\infty\bigl[W_\ell(r)f(r)^2+(W_\ell(r)+E_\ell(r))k(r)^2\bigr]r\,dr\\
	&\quad\quad\leq\liminf_{j\to\infty}\int_0^\infty\bigl[W_\ell(r)f_j(r)^2+(W_\ell(r)+E_\ell(r))k_j(r)^2\bigr]r\,dr\\
	&\quad\quad\leq\liminf_{j\to\infty}\left(\mathscr R_\ell(f_j)+\mathscr R_\ell(k_j)+\mathscr Q_\ell(k_j)
	-4\ell\int_0^\infty\frac{B(r)}r f_j(r)k_j(r)\,dr\right)\\
	&\quad\quad=\mathscr R_\ell(f)+\mathscr R_\ell(k)+\mathscr Q_\ell(k)
	-4\ell\int_0^\infty\frac{B(r)}r f(r)k(r)\,dr.
\end{aligned}
$$
Here, we have used the fact that $B(r)/r^2$ is bounded in $(0,\infty)$ to derive the last equality. Since $W_\ell(r)>0$ and $W_\ell(r)+E_\ell(r)>0$ for $r\neq1$, we know \eqref{inequality higher l} holds for $\ell=1,2$, and the equality holds if and only if $f=k=0$ almost everywhere.

\subsection{The case $\ell\geq3$}\label{subsection 4.4}
For $\ell\geq 3$, define the following coefficient matrix
\begin{equation}\label{coefficient matrix}
	V_\ell(r)\coloneqq
	\begin{pmatrix}
		r^2V(r)+(1-B(r))\ell^2&-2\ell B(r)\\
		-2\ell B(r)&r^2(V(r)-\alpha^2)+(2-B(r))\ell^2
	\end{pmatrix}.
\end{equation}
Then the left-hand side of \eqref{inequality higher l} can be written as
\begin{equation}\label{Dl another form}
	\int_0^\infty\bigl((1-B)f'^2+(2-B)k'^2\bigr)r\,dr+\int_0^\infty\begin{pmatrix}f&k\end{pmatrix} V_\ell(r)\begin{pmatrix}f\\k\end{pmatrix}\frac{dr}{r}.
\end{equation}
The following lemma is crucial.
\begin{lemma}\label{positive definite}
	For $\ell\geq3$, the matrix $V_\ell(r)$ is positive definite for any $r\in(0,1)\cup(1,\infty)$.
\end{lemma}
\begin{proof}
	First, we consider the inner part $0<r<1$. Write $s=r^2$. By \eqref{log derivative identity} and \eqref{definition of V}, we may write
	$$
	V(\sqrt s)=-\alpha^2p(s)+2sC(s)\bigl[2+a(s)(1-s)\bigr].
	$$
	Denote
	$$
		A_3(s)\coloneqq(V_3(\sqrt s))_{11}=sV(\sqrt s)+9p(s)
	$$
	and
	$$
	\Delta(s)\coloneqq(V_3(\sqrt s))_{22}-(V_3(\sqrt s))_{11}=9-\alpha^2s.
	$$
	Then for any $\ell\geq3$, it follows from \eqref{coefficient matrix} that
	$$
	(V_\ell(\sqrt s))_{11}=A_3(s)+p(s)(\ell^2-9)
	$$
	and
	$$
	(V_\ell(\sqrt s))_{22}=A_3(s)+\Delta(s)+(1+p(s))(\ell^2-9).
	$$
	Since $\alpha^2<29/5$, we know
	$$
	\Delta(s)\geq\Delta_0(s)\coloneqq 9-\frac{29}{5}s>0.
	$$
	By Lemma \ref{lemma on delta and a} and \ref{estimate of c(s) and p(s)}, we have
	$$
	\begin{aligned}
		A_3(s)&=\Delta(s)p(s)+2s^2C(s)\bigl[2+a(s)(1-s)\bigr]\\
		&\geq A_0(s)\coloneqq\Delta_0(s)p_0(s)+s^2\left[4+\frac{89}{50}(1-s)\right]\left[\frac{13}{10}+\frac{17}{20}(1-s)\right],
	\end{aligned}
	$$
	where $p_0(s)$ is given by \eqref{estimate of p(s) equation}. We will use the polynomials $\Delta_0(s)$ and $A_0(s)$ to bound $\Delta(s)$ and $A_3(s)$ from below. Specifically, we claim that
	\begin{equation}\label{diagonal element estimate}
		A_0(s)+\Delta_0(s)\ge\frac{42}{5},\quad A_0(s)>\frac{30}{7}\bigl(1-p_0(s)\bigr)^2,\quad 0\leq s\leq1.
	\end{equation}
	By a direct expansion, we have
	$$
		A_0(s)=\frac{26}{5}-\frac{2343}{500}(1-s)-\frac{527}{200}(1-s)^2+\frac{11047}{1500}(1-s)^3+\frac{9469}{3000}(1-s)^4.
	$$
	Hence, 
	$$
	A_0(s)+\Delta_0(s)-42/5\ge(1-s)\bigl[1-3(1-s)+7(1-s)^2\bigr]\geq0.
	$$
	On the other hand, since $0\leq 1-s\leq 1$, we know
	$$
	0\leq p_0(s)\leq \frac{14}{15}(1-s)^2,
	$$
	and hence
	$$
	\begin{aligned}
		A_0(s)-\frac{30}{7}\bigl(1-p_0(s)\bigr)^2&\geq\frac{32}{35}-\frac{2343}{500}(1-s)+\frac{4111}{1400}(1-s)^2\\
		&\quad+(1-s)^3\left[\frac{102829}{10500}-\frac{577}{1000}(1-s)\right]\\
		&\geq\frac{32}{35}-\frac{2343}{500}(1-s)+\frac{4111}{1400}(1-s)^2+\frac{193541}{21000}(1-s)^3.
	\end{aligned}
	$$
	Since
	$$
	(1-s)^3-\frac{1-s}{3}+\frac{2}{27}=\left(s-\frac23\right)^2\left(\frac53-s\right)\geq0,
	$$
	we obtain that
	$$
	\begin{aligned}
		A_0(s)-\frac{30}{7}\bigl(1-p_0(s)\bigr)^2&\geq\frac{23-162(1-s)+293(1-s)^2}{100}\\
		&=\frac{(293s-212)^2+178}{29300}>0.
	\end{aligned}
	$$
	This gives \eqref{diagonal element estimate}.
	
	Since $0\leq b(s)\leq 1-p_0(s)$, the diagonal entries satisfy
	$$
	(V_\ell(\sqrt s))_{11}\geq A_0(s)>\frac{30}{7}b(s)^2
	$$
	and
	$$
	(V_\ell(\sqrt s))_{22}\geq A_0(s)+\Delta_0(s)+(\ell^2-9)\geq\ell^2-\frac35\ge\frac{14}{15}\ell^2.
	$$
	Therefore, we know
	$$
	(V_\ell(\sqrt s))_{11}(V_\ell(\sqrt s))_{22}>4\ell^2b(s)^2=(V_\ell(\sqrt s))_{12}(V_\ell(\sqrt s))_{21}.
	$$
	More precisely, from the above estimates, we have
	$$
	\det V_\ell(r)\geq\frac{14\ell^2}{15}\cdot \frac{(293s-212)^2+178}{29300}>0.
	$$
	Clearly, for $0<s<1$,
	$$
		(V_\ell(\sqrt s))_{11}>\frac{30}{7}b(s)^2>0,
	$$
	so $V_\ell(r)=V_\ell(\sqrt s)$ is positive definite in $(0,1)$.
	
	For $r>1$, we know
	$$
	V_\ell(r)=
	\begin{pmatrix}
		dr^2-1&-2\ell\\
		-2\ell&(\beta^2-2\alpha^2)r^2+\ell^2-1
	\end{pmatrix}.
	$$
	Both diagonal entries are positive and strictly increasing with respect to $r>1$. By a direct calculation, we know for any $r>1$, there holds
	$$
	\det V_\ell(r)>\lim_{t\downarrow 1}\det V_\ell(t)=(d-1)(\beta^2-2\alpha^2-1)+(d-5)\ell^2>0,
	$$
	from which we know $V_\ell(r)$ is positive definite.
\end{proof}

By \eqref{Dl another form} and Lemma \ref{positive definite}, we know that \eqref{inequality higher l} holds for every $\ell \geq 3$, with equality holding if and only if $f=k=0$ almost everywhere.

\subsection{Completion of the proof: summation and equality}
Finally in this subsection, we give the proof of Proposition \ref{two function proposition}.
\begin{proof}[Proof of Proposition \ref{two function proposition}]
	Write the angular Fourier expansions
	$$
	\begin{aligned}
		u(r,\theta)	&=u_0(r)+\sum_{\ell=1}^{\infty}u_\ell(r,\theta)=u_0(r)+\sum_{\ell=1}^{\infty}\bigl(u_\ell^c(r)\cos(\ell\theta)+u_\ell^s(r)\sin(\ell\theta)\bigr),\\
		v(r,\theta)&=v_0(r)+\sum_{\ell=1}^{\infty}v_\ell(r,\theta)=v_0(r)+\sum_{\ell=1}^{\infty}\bigl(v_\ell^c(r)\cos(\ell\theta)+v_\ell^s(r)\sin(\ell\theta)\bigr).
	\end{aligned}
	$$
	Recall that we have proved in Section \ref{subsection 4.2} that $D_0\geq 0$. Moreover, by taking $(f,k)=(u_\ell^c,v_\ell^s)$ and $(f,k)=(u_\ell^s,-v_\ell^c)$, we have proved that $D_\ell\geq 0$ for all $\ell\geq1$. 
	
	Since $S_Nu\to u$ and $S_Nv\to v$ in $H^1(\R^2)$, we know from the finite Fourier decomposition and the continuity that
	$$
	\begin{aligned}
		&R_0(u)+R_0(v)+q_{\alpha^2}^{\R^2}(v)-\mathcal{J}_0(u,v)\\
		&\quad\quad=\lim_{N\to\infty}\bigl[R_0(S_Nu)+R_0(S_Nv)+q_{\alpha^2}^{\R^2}(S_Nv)-\mathcal{J}_0(S_Nu,S_Nv)\bigr]\\
		&\quad\quad=\lim_{N\to\infty}\sum_{\ell=0}^{N}D_\ell =\sum_{\ell=0}^{\infty}D_\ell\geq0.
	\end{aligned}
	$$
	The equality holds if and only if $D_\ell=0$ for all $\ell\geq0$. By the analysis presented earlier in this section, we know
	$$
	D_0=0 \iff (u_0,v_0)=(c\phi,0)\quad\text{for some constant }c\in\R
	$$
	and
	$$
	D_\ell=0\iff (u_\ell,v_\ell)=(0,0),\quad\ell\geq1.
	$$
	Hence the equality in \eqref{two function equation} holds if and only if $(u,v)=(c\phi,0)$ for some $c\in\R$. The proof is complete.
\end{proof}

\section*{Acknowledgements} Zhen Song is supported by LiaoNing Revitalization Talents Program (No. XLYC2502044). Wenming Zou is supported by National Key R\&D Program of China (Grant 2023YFA1010001) and National Natural Science Foundation of China (No. 12571123 and 12671139).

\section*{Data Availability} No data was used for the research described in this paper.

\section*{Statements and Declarations}

\subsection*{Conflict of interest} On behalf of all authors, the corresponding author states that there is no conflict of interest.

\subsection*{AI usage statement} During the preparation of this work, the authors used ChatGPT models as interactive research-assistance tools for literature reading, exploring proof strategies, calculations, checking arguments, and manuscript preparation. The paper presented here is the fruit of successive revisions of the original one, and has been written entirely by the authors. We take full responsibility for all mathematical statements, arguments, citations, and conclusions in this paper.

	\bibliographystyle{amsplain}

\end{document}